\documentclass[leqno, 12pt]{article} 

\usepackage[latin1]{inputenc} 
\usepackage{amsmath,amssymb} 
\usepackage{latexsym} 

\newtheorem{theorem}{Theorem}[section] 
\newtheorem{corollary}[theorem]{Corollary} 
\newtheorem{lemma}[theorem]{Lemma} 
\newtheorem{proposition}[theorem]{Proposition}

\newtheorem{remark}[theorem]{Remark} 
\numberwithin{equation}{section}

\def \ov{\overline} 
 
\def \be{\begin{equation}} 
\def \ee{\end{equation}} 
\def \bt{\begin{theorem}} 
\def \et{\end{theorem}} 
\def \bea{\begin{eqnarray}} 
\def \eea{\end{eqnarray}} 
\def \bas{\begin{eqnarray*}} 
\def \eas{\end{eqnarray*}} 
 
\def \bb{\beta} 
\def \ga{\gamma} 
 
\def \de{\delta}

\def \la{\lambda}

\def \th{\theta}

\def \ff{\infty} 
\def \wh{\widehat} 
\def \wt{\widetilde} 
 
\def \rar{\rightarrow}

\def \R{{\bf R}}

\def \Z{{\bf Z}}

\def \HH{{\cal H}}

\def \({\left(} 
\def \){\right)}

\def \lc{\left\{} 
\def \rc{\right\}}

\def \nn{\nonumber}

\def \bc{\begin{center} } 
\def \ec{\end{center} } 
\def \bs{\begin{slide} } 
\def \es{\end{slide} } 
 
\def\R{{\mathbb R}} 
\def\Z{{\mathbb Z}}

\def\eps{\varepsilon}

\def\lam{\lambda}

\def\phi{\varphi} 
\def\proof{{\medskip\noindent {\bf Proof. }}} 
 
\def\qed{{\hfill $\square$ \bigskip}}

\def\square{{\vcenter{\vbox{\hrule height.3pt 
         \hbox{\vrule width.3pt height5pt \kern5pt 
            \vrule width.3pt} 
         \hrule height.3pt}}}} 
 
\def\sA {{\cal A}}

 \def\bE {{\Bbb E}}

\def\bP {{\Bbb P}}

\def\square{{\vcenter{\vbox{\hrule height.3pt 
          \hbox{\vrule width.3pt height5pt \kern5pt 
             \vrule width.3pt} 
          \hrule height.3pt}}}} 
\def\qed{{\hfill $\square$ \bigskip}} 
\def\ol{\overline} 
\def\E{{\bE}} 
\def\P{{\bP}}

\begin{document} 
 
\title{Moderate deviations for the range 
of planar 
   random walks} 
\author{Richard Bass\thanks{Research partially supported by NSF grant 
\#DMS-0244737} 
    \, \, Xia Chen\thanks {Research   partially supported by NSF grant 
\#DMS-0405188.}\, \, Jay Rosen\thanks {Research partially supported 
by grants from the NSF  and from PSC-CUNY.}} 
 
%\date{October 6, 2004} 
 
\maketitle 
 
\bibliographystyle{amsplain} 
 
\begin{abstract} Given a symmetric random walk in $\Z^2$ with finite 
second moments, let $R_n$ be the range of the random walk up to time 
$n$.  We study moderate deviations for $R_n -\E R_n$ and $\E R_n -R_n$. 
We also derive the corresponding laws of the iterated logarithm. 
\end{abstract}

\section{Introduction} 
 
Let $X_i$ be symmetric i.i.d.\ random vectors taking values in $\Z^2$ 
with mean 0 and finite covariance matrix $\Gamma$, set 
$S_n=\sum_{ i=1}^{ n}X_i$, and suppose that no proper subgroup of
$\Z^2$ supports the random walk $S_n$.
For any 
random variable $Y$ we will use the notation 
\[ 
\ol Y=Y -\E Y. 
\] 
\medskip Let 
\begin{equation} R_n=\#\{S_1, \ldots, S_n\} 
\label{defrange} 
\end{equation} be the range of the random walk up to time $n$. The 
purpose of this paper is to obtain moderate deviation results for $\ol 
R_n$ and $-\ol R_n$. 
 
\medskip For moderate deviations of $\ol R_n$ we have the following. 
Let
\begin{equation}\label{Hdef}
\HH(n)=\sum_{k=0}^n \P^0(S_k=0).
\end{equation}
Since the $X_i$ have two moments, then by \cite{LR}, Section 2,
$$ \HH(n)=\sum_{k=0}^n \P^0(S_k=0) \sim \frac{\log n}{2\pi \sqrt{\det \Gamma}}$$
and
$$\HH(n)-\HH([n/b_n])=\sum_{k=[n/b_n]+1}^n \P^0(S_k=0)  
\sim \frac{\log b_n}{2\pi \sqrt{\det \Gamma}}.$$
 
\begin{theorem}\label{theo-range} Let $\{b_n\}$ be a positive sequence 
satisfying $b_n\to \infty$ and  
$\log b_n=o((\log n)^{1/2})$ as $n\to \infty$.  
There are two constants $C_1, C_2>0$ independent of the choice of the 
sequence 
$\{b_n\}$ such that 
\begin{align} -C_1&\le\liminf_{n\to\infty}b_n^{-1}\log 
\P\Big\{\ol{R}_n 
\ge {n\over\HH( n)^2}(\HH(n)-\HH([n/b_n]))\Big\}\nn\\ 
&\le \limsup_{n\to\infty}b_n^{-1}\log \P\Big\{\ol{R}_n 
\ge {n\over\HH( n)^2}(\HH(n)-\HH([n/b_n])\Big\}\le 
-C_2. 
\label{erw21} 
\end{align} 
\end{theorem}

\begin{remark}\label{remrange}
{\rm The proof will show that $C_2$ in the statement of Theorem 
\ref{theo-range} is equal to the constant 
$L$ given in Theorem 1.3 in \cite{BC}. We believe that $C_1$ is also
equal to $L$, but we do not  have a proof of this fact. 
}
\end{remark}

A more precise statement than Theorem \ref{theo-range} is possible when
the $X_i$ have slightly more than two moments.

\begin{corollary}\label{corrange}
Suppose $\E[\,|X_i|^2 ( \log^+ (|X_i|))^{\frac12+\delta}\,]<\infty$
for some $\delta>0$.
Let $\{b_n\}$ be a positive sequence
satisfying $b_n\to \infty$ and
$\log b_n=o((\log n)^{1/2})$ as $n\to \infty$.
There are two constants $C_1, C_2>0$ independent of the choice of the
sequence
$\{b_n\}$ such that
\begin{align} -C_1&\le\liminf_{n\to\infty}b_n^{-\theta}\log
\P\Big\{\ol{R}_n
\ge 2\theta \pi\sqrt{\det\Gamma}{n\over(\log n)^2}\log b_n\Big\}\nn\\
&\le \limsup_{n\to\infty}b_n^{-\theta}\log \P\Big\{\ol{R}_n
\ge 2\theta\pi\sqrt{\det\Gamma}{n\over(\log n)^2}\log b_n\Big\}\le
-C_2
\label{erw21B}
\end{align} for any $\theta>0$.
\end{corollary}

\begin{remark}\label{remrange2}
{\rm The constants $C_1$, $C_2$ are the same as in the statement of  Theorem
\ref{theo-range}. See Remark \ref{remrange}.
}
\end{remark}

For $b_n$ tending to infinity faster than the rate given in
Theorem \ref{theo-range}, e.g., $\log b_n=(\log n)^2$,  
then  
we are in the realm of large 
deviations. For results on large deviations of the range, see 
\cite{DV}, \cite{HamanaKesten1}, \cite{HamanaKesten2}.

For the moderate deviations of $-\ol R_n=\E R_n-R_n$ we have the following. Let 
$\kappa(2,2)$ be the smallest $A$ such that 
$$\| f\|_4\leq A \| \nabla f\|_2^{1/2} \|f\|_2^{1/2}$$ for all $f\in C^1$ 
with compact support. (This constant appeared in \cite{BC}.) 
 
\begin{theorem}\label{LU1.3} Suppose  $b_n\to \infty$ and $b_n 
=o((\log n)^{1/5})$ as $n\to \infty$. For $\lambda>0$ 
$$\lim_{n\to \infty} \frac{1}{b_n} \log 
\P\Big(-\ol R_n>\lambda \frac{nb_n}{\log^2 n}\Big) =-(2\pi)^{-2} (\det 
\Gamma)^{-1/2} \kappa(2,2)^{_4}\lambda.$$ 
\end{theorem}

Comparing Theorems \ref{theo-range} and \ref{LU1.3}, we see that  the 
upper and lower tails of $\ol R_n$ are quite different. This is similar to 
the behavior of the distribution of the self-intersection local time of 
planar Brownian motion. This is not surprising, since LeGall,  
\cite[Theorem 6.1]{LeGall6}, shows  that $\ol R_n$, properly normalized, converges in 
distribution to the self-intersection local time. 
 
The moderate deviations of $\ol R_n$ are quite similar in nature to those 
of $-\ol L_n$, where $L_n$ is the number of self-intersections of the 
random walk $S_n$; see \cite{BCR05}. Again, \cite[Theorem 6.1]{LeGall6} gives a partial 
explanation of this. However the case of the range is much more 
difficult than the  corresponding results for intersection local times. The 
latter case can be represented as a quadratic functional of the path, which is 
amenable to the techniques of large deviation theory, while the range 
cannot be so represented. This has necessitated the development of several new tools, see in particular Sections \ref{sec-smoothrange} and \ref{sec-expapp}, which we expect will have further applications in the study of the range of random walks.

\medskip Theorem  \ref{theo-range}  gives rise to the following LIL for 
$\ol R_{  n}$. 
 
\begin{theorem}\label{theo-lilrw} 
\begin{equation} 
\limsup_{n\to\infty}\frac{\ol R_n}{n\log \log\log n/\log^2 n} = 2\pi 
\sqrt{\det \Gamma}, \qquad 
\mbox{\rm  a.s.} 
\label{19.53rw} 
\end{equation} 
\end{theorem} 
 
This result is an improvement of that in  \cite{BK}; there it was required 
that the 
$X_i$ be bounded random variables and the constant was not identified. 
Theorem \ref{theo-range} is a more precise estimate than is needed for
Theorem \ref{theo-lilrw}; this is why Theorem \ref{theo-range} needs to
be stated in terms of $\HH(n)$ while Theorem \ref{theo-lilrw} does not.
 
For an LIL for $-\ol R_n$ we have a different rate. 
 
\begin{theorem}\label{LT5.3} We have 
$$\limsup_{n\to \infty} \frac{-\ol R_n}{n\log\log n/\log^2 
n}=(2\pi)^{-2} \sqrt{\det \Gamma}\,\, \kappa(2,2)^4, 
\qquad \hbox{\rm a.s.}$$ 
\end{theorem}

The study of the range of a lattice-valued (or $\Z^d$-valued) 
random walk has a long history in probability and the results show 
a strong dependence on the dimension $d$. 
See \cite{DvEr}, \cite{JaPr1}, \cite{LeGall6}, \cite{LR}, 
\cite{HamanaKesten1}, \cite{HamanaKesten2}, 
\cite{DV}, and \cite{BK} and the references in these 
papers, to cite only a few.  The two dimensional case seems to be the most 
difficult; in one dimension no renormalization is needed (see \cite{Cc}), 
while for $d\geq 3$ the tails are sub-Gaussian and have asymptotically 
symmetric behavior. In two dimensions, renormalization is needed 
and the tails have non-symmetric behavior. In this case,  
the central limit theorem was proved in 1986 in \cite{LeGall6}, 
while the first law of the iterated logarithm was not proved until a few 
years ago in \cite{BK}.

\medskip
\noindent{\bf Acknowledgment:} We would like to thank Greg Lawler
and Takashi Kumagai for helpful discussions and their interest in
this paper.

\section{Moments of the range}

In this section we first give an estimate for the expectation of
the range.

By \cite{LR}, Theorem 6.9, we have
\begin{equation}\label{lrth6.9}
\E R_n= \frac{n}{\HH(n)} +\frac{1}{2\pi \sqrt{\det \Gamma}} \frac{n}{\HH(n)^2}(1+o(1)),
\end{equation}
where $\HH$ is defined in (\ref{Hdef}).
By \cite{LR}, Section 2, 
\begin{equation} \label{secr11}
\HH(n)\sim \frac{\log n}{2\pi \sqrt{\det \Gamma}}
\end{equation}
and 
\begin{equation} \label{secr12}
\HH(n)-\HH(m)\sim \frac{\log(n/m)}{2\pi \sqrt {\det \Gamma}}
\end{equation}
as $n$ and $m$ tend to infinity.

Throughout this paper we will mostly be concerned with random walks that
have only second moments. The exception is the following proposition,
which  supposes slightly more than two  moments, and Corollary \ref{corrange}.

\begin{proposition}\label{p21a} Suppose $\{X_i\}$ is a sequence of i.i.d.\
mean  zero random vectors taking values in $\Z^2$ with
\begin{equation}
\E\(|X|^{ 2}(\log^{ +} |X|)^{ \frac12+\de}\)<\ff\label{2.1}
\end{equation}
for some $\delta>0$ and nondegenerate
covariance matrix $\Gamma$. Let
$S_n=\sum_{i=1}^n X_i$ and suppose $S_n$ is strongly aperiodic. Then
\begin{equation}\label{2.2}
\P(S_n=0)={ 1\over 2\pi n\sqrt{\det
\Gamma}} +O\Big({1 \over n(\log n)^{(1+\de)/2}}\Big).
\end{equation}
\end{proposition}

\proof  Let $\phi$ be the characteristic function of $X_i$, let 
$x\cdot y$ denote
the inner product in $\R^2$, let
$Q(u)=u\cdot \Gamma u$, and let $C=[-\pi,\pi]^2$.   We observe that
\begin{align} 
|1-\phi( u)&-Q(u)|\label{2.3}\\ &
=|\E\(1-e^{ iu\cdot X} +iu\cdot X+( 1/2)( iu\cdot X)^{ 2}\)| \nonumber\\ &
\leq c_1 |u|^{ 3}\E\(1_{ \{|X|\leq 1/|u| \}}|X|^{ 3}\) + c_1 |u|^{ 
2}\E\(1_{ \{|X|> 1/|u|
\}}|X|^{ 2}\) \nonumber
\end{align}
and consequently for any fixed $M>0$
\begin{align} \qquad
|1&-\phi( u/\sqrt n)-Q(u/\sqrt n)|\label{2.4}\\  &
\leq c_2 \({1\over n^{ 3/2}}\)\E\(1_{ \{|u| |X|\leq \sqrt n\}} (|u||X|)^{ 3}\)
+ c_2\({1\over n}\)\E\(1_{\{|u| |X|>\sqrt n\}}(|u||X|)^{ 2}\)\nonumber\\  &
\leq c_3{1\over n^{ 3/2}}+c_3 \({1\over n^{ 3/2}}\)\E\(1_{ \{M<|u| |X|\leq \sqrt
n\}} (|u||X|)^{ 3}\)
  \nonumber\\  &\hspace{ 1in}
+ c_3\({1\over n}\)\E\(1_{\{|u| |X|>\sqrt n\}}(|u||X|)^{ 2}\). \nonumber
\end{align}
Choose M so that $x/\log^{ 1/2+\de}(x)$ is monotone increasing on $x\geq
M$, and therefore
\begin{align} 
\E(&1_{ \{M<|u| |X|\leq \sqrt n\}} (|u||X|)^{ 3})\label{2.5}\\ &
\leq \E\(1_{ \{M<|u| |X|\leq \sqrt n\}} (|u||X|)^{ 2}\log^{ 1/2+\de}(|u||X|)
{ |u||X|\over \log^{1/2+\de}(|u||X|)}\)\nonumber\\ &
\leq \({ \sqrt n\over \log^{ 1/2+\de}(\sqrt n)}\)\E\(1_{ \{M<|u| |X|\leq \sqrt
n\}} (|u||X|)^{ 2}\log^{ 1/2+\de}(|u||X|) \).\nonumber
\end{align}
Also
\begin{align} 
\E(&1_{\{|u| |X|>\sqrt n\}}(|u||X|)^{ 2})\label{2.6}\\ &
\leq \({1\over \log^{ 1/2+\de}(\sqrt n)}\)\E\(1_{\{|u| |X|>\sqrt n\}}(|u||X|)^{
2}\log^{ 1/2+\de}(|u||X|)\).\nonumber
\end{align}
(\ref{2.4}) then implies that
\begin{eqnarray} &&
|1-\phi( u/\sqrt n)-Q(u/\sqrt n)|
\leq c { |u|^{ 2}|\log^{ 1/2+\de}(|u|)|\over n\log^{1/2+\de}( n)}.
\label{2.7}
\end{eqnarray}

  Following the proof in Spitzer \cite{spitzer}, pp.~76--77,
\begin{align*} 2\pi n \P(S_n=0)&=(2\pi)^{-1} \int_{\sqrt n C} \phi(u/\sqrt
n)^n du\\ &=I_0+I_1(n,A_n)+I_2(n,A_n)+I_3(n,A_n,r)+I_4(n,r),
\end{align*} where
\begin{align*} I_0&=(2\pi)^{-1}\int_{\R^2} e^{-Q(u)/2} du=(\det Q)^{-1/2},\\
I_1(n,A_n)&=(2\pi)^{-1}\int_{|u|\leq A_n} [\phi(u/\sqrt n)^n -e^{-Q(u)/2}]\,
du,\\ I_2(n,A_n)&=-(2\pi)^{-1}\int_{|u|>A_n} e^{-Q(u)/2} du,\\
I_3(n,A_n,r)&=(2\pi)^{-1} \int_{A_n<|u|<r\sqrt n} \phi(u/\sqrt n)^n\, du,\\
I_4(n,r)&=(2\pi)^{-1} \int_{|u|\geq r\sqrt n, u\in \sqrt n C}
\phi(u/\sqrt n)^n\, du.
\end{align*} Still following \cite{spitzer}, we can choose $r$ such that
$|\phi(u/\sqrt n)^n|\leq e^{-Q(u)/4}$ if
$|u|\leq r\sqrt n$ and by the strong aperiodicity there exists $\gamma>0$ such
that
$|\phi(u/\sqrt n)|\leq 1-\gamma$ if $|u|>r\sqrt n$ and $u\in \sqrt n C$. Set
$A_n=c_4\sqrt{\log \log n}$. We have
$$|I_4(n,r)|\leq (2\pi)^{-1}\int_{u\in \sqrt n C} (1-\gamma)^n du=O(n^{-p})$$
for every positive integer $p$. Next
$$|I_3(n,A_n,r)|\leq \int_{|u|>c_4\sqrt{\log \log  n}} e^{-Q(u)/4}du=O((
\log n)^{-2})$$ for
  $c_4$ large and similarly we have the same bound for
$|I_2(n,A_n)|$. To estimate $I_1(n,A_n)$ we use the  inequality $|a^n-b^n|\leq
n|a-b|$ if $|a|, |b|\leq 1$ with $a=\phi(u/\sqrt n)$ and 
$b=e^{-Q(u)/2n}$. Using
(\ref{2.7}) and the analogous expansion for $e^{-Q(u)/2n}$ we
have
\begin{align*} |\phi(u/\sqrt n)^n- e^{-Q(u)/2}|&\leq n|\phi(u/\sqrt
n)-e^{-Q(u)/2n}|\\ &\leq c_5n { |u|^{ 2}|\log^{ 1/2+\de}(|u|)|\over n\log^{
1/2+\de}( n)} =c_5 { |u|^{ 2}|\log^{ 1/2+\de}(|u|)|\over \log^{
1/2+\de}( n)}.
\end{align*} Integrating this over the set $\{|u|\leq A_n\}$, we see
$$|I_1(n,A_n)|=O((\log \log  n)^{2+\delta/2}/(\log
n)^{1/2+\de})=O(1/(\log n)^{(1+\de)/2}).$$ Summing $I_0$ through $I_4$,
we obtain
$$2\pi n\P(S_n=0)=(\det \Gamma)^{-1/2} + O(1/(\log  n)^{(1+\de)/2}).$$
\qed

\def\be{{\bf e}}

Next we 
establish some sharp exponential 
estimates for the range and intersection of ranges. Aside from their 
intrinsic interest,  they will be used to estimate the tail probabilities in 
our first main theorem.   

We write $S(I)$ for $\{S_k: k\in I\}$. 
Let $S^{(i)}$, $i=1, \ldots, p$ be 
$p$ independent copies of $S$. First, by  Corollary 1 of 
\cite{Cb}, for any integers $a\ge 1$, $n_1,\cdots n_a\ge 1$, 
\begin{equation} 
\big(\E J^m_{n_1+\cdots +n_a}\big)^{1/p}\le 
\sum_{\scriptstyle k_1+\cdots +k_a =m\atop\scriptstyle k_1,\cdots, 
k_a\ge 0}{m!\over k_1!\cdots k_a!} 
\big(\E J_{n_1}^{k_1}\big)^{1/p}\cdots \big(\E 
J_{n_a}^{k_a}\big)^{1/p},\label{erw14bb} 
\end{equation} where 
$$ J_n=\#\big\{S^{(1)}[1,n]\cap\cdots \cap S^{(p)}[1,n]\big\} 
\hskip.2in n=1,2,\cdots. 
$$

In the next Theorem we deduce from this the  exponential integrability of $J_n$, which was 
established in \cite{BK} in the special case $p=2$ and under the 
condition that $S$ had bounded increments. 
 
\begin{theorem}\label{rd12} Assume that the planar random walk $S$ has finite 
second moments and zero mean. There exists $\theta >0$ such that 
\begin{equation} 
\sup_{n}\sup_{y_1,\cdots,y_p}\E^{(y_1,\cdots, y_p)}\exp\Big\{\theta 
\Big({(\log n)^p\over n}\Big)^{1/(p-1)}J_n^{1/(p-1)}\Big\} 
<\infty.\label{erw19} 
\end{equation} 
\end{theorem}

\proof  
We recall  the fact (see Remarks, p.~664, in \cite{LR}) that 
\begin{equation} 
\E J_n^{k}\le (k!)^p(\E J_n)^k,\hskip.2in k=0,1,\cdots,\label{erw14a} 
\end{equation} and for some $C<\infty$ 
\begin{equation} 
\E J_n\le  {C n\over (\log n)^p},\hskip.2in n=1,\cdots.\label{erw14b} 
\end{equation} The proof of (\ref{erw19}) is a modification of the 
approach used in  Lemma 1 of \cite{Cb}. We begin by  showing that there is 
a constant $C>0$ such that 
\begin{equation} 
\sup_{n} \E J_n^m \le C^m(m!)^{p-1}\Big({n\over(\log n)^p}\Big)^m, \hskip.2in 
m,n=1,2,\cdots.\label{erw111} 
\end{equation} We first consider the case $m\le (\log n)^{(p-1)/p}$. 
Write 
   $l(n,m)=[n/m]+1$. Then by (\ref{erw14bb}) and (\ref{erw14b}),  
\begin{align*} 
\big(\E J_n^m\big)^{1/p}&\le 
\sum_{\scriptstyle k_1+\cdots +k_m =m\atop\scriptstyle k_1,\cdots, 
k_m\ge 0}{m!\over k_1!\cdots k_m!} 
\big(\E J_{l(n,m)}^{k_1}\big)^{1/p}\cdots 
\big(\E J_{l(n,m)}^{k_m}\big)^{1/p}\\ &\le \sum_{\scriptstyle 
k_1+\cdots +k_m =m\atop\scriptstyle k_1,\cdots, k_m\ge 0}{m!\over 
k_1!\cdots k_m!}k_1!\cdots k_m! 
\big(\E J_{l(n,m)}\big)^{k_1/p}\cdots 
\big(\E J_{l(n,m)}\big)^{k_m/p}\\ &={2m-1 \choose m}   m! 
\Big(\E J_{l(n,m)}\Big)^{m/p} 
\le {2m-1 \choose m} m! C^m \Big({(n/m)\over (\log n)^p}\Big)^{m/p}\\ 
&\le {2m \choose m}  (m!)^{p-1\over p}C^m 
\Big({n\over (\log n)^p}\Big)^{m/p}, 
\end{align*} where the second inequality follows from (\ref{erw14a}) and 
the third from (\ref{erw14b}) using  
the fact that 
$m=O( \log n)$ so that $\log n=O( \log (n/m))$. Hence, taking $p$-th powers  we obtain 
$$ 
\E J_n^m\le {2m \choose m}^pC^{pm}(m!)^{p-1} 
\Big({n\over (\log n)^p}\Big)^m, 
$$ and (\ref{erw111}) for the case of $m\le (\log n)^{(p-1)/p}$ follows 
from the fact 
$$ {2m \choose m}\le 4^m. 
$$ 
 
For the case $m>(\log n)^{(p-1)/p}$, notice from the definition of $J_n$ 
   that 
$J_n\le n$. So we have 
\begin{align*} 
\E J_n^m\le n^m &=(\log 
n)^{pm}\Big({n\over (\log n)^p}\Big)^m 
\le m^{(p-1)m}\Big({n\over (\log n)^p}\Big)^m\\ &\le (m!)^{p-1}C^m 
\Big({n\over (\log n)^p}\Big)^m, 
\end{align*} 
  where the last step follows from Stirling's formula. This completes the proof of (\ref{erw111}). 
\medskip
 
By H\"{o}lder's inequality this shows that
\begin{eqnarray}
&& 
\Big({(\log n)^p\over n}\Big)^{m/(p-1)}\sup_{y_1,\cdots,y_p}\E^{(y_1,\cdots, y_p)} \Big(J_n^{m/(p-1)}\Big)
\label{ 2.6add }\\
&&   
\leq \Big({(\log n)^p\over n}\Big)^{m/(p-1)}\sup_{y_1,\cdots,y_p}\lc \E^{(y_1,\cdots, y_p)} \Big(J_n^{m}\Big)\rc^{1/(p-1)}\nonumber\\
&&   
\leq \Big({(\log n)^p\over n}\Big)^{m/(p-1)}\lc \E \Big(J_n^{m}\Big)\rc^{1/(p-1)}
\leq C^{m}m!\nonumber
\end{eqnarray} 
where the second inequality  used  \cite{Cb}, p.1053. Our theorem then follows from a Taylor expansion. 
\qed 
 
{\bf Remark.} Theorem \ref{rd12} is sharp in the sense that 
(\ref{erw19}) does not hold if $\theta$ is too large. 
Indeed, by \cite{LeGall6}, 
for any $m=1,2,\cdots$, 
$$ {(\log n)^{pm}\over n^m}\E J_n^m\longrightarrow (2\pi)^{pm}\det 
(\Gamma)^{m/2} 
\E\alpha\big([0,1]^p\big)^m 
$$ as $n\to \infty$, where $\alpha \big([0,1]^p\big)$ is the Brownian 
intersection local time formally defined by 
$$ 
\alpha \big([0,1]^p\big)=\int_{\R^d}\bigg[\prod_{j=1}^p\int_0^1 
\delta_x\big(W_j(s)\big)ds\bigg]dx, 
$$ and by Theorem 2.1 in \cite{C} 
$$ 
\E\exp\Big\{\theta\alpha \big([0,1]^p\big)^{(p-1)^{-1}} \Big\}=\infty 
$$ for large $\theta $. The following theorem  is sharp in the same sense.

\begin{theorem} \label{rd13} Assume that the planar random walk $S$ 
has finite second moments and zero mean. Then there exists $\theta >0$ 
such that 
\begin{equation} 
\sup_n\E\exp\Big\{\theta{(\log n)^2\over n}\vert\ol{R}_n\vert\Big\} 
<\infty.\label{erw112} 
\end{equation} 
\end{theorem} 
 
\proof We first consider the case where $n$ is replaced by $2^n$. Let 
$$ N= [2(\log 2)^{-1}\log n] 
$$ so that $2^N\sim n^{ 2}$ and note that
\begin{eqnarray}
&&\hspace{.3 in}\#\big\{S[1,2^n]\big\} =\sum_{k=1}^{2^N}\#\big\{S((k-1)2^{n-N}, 
k2^{n-N}]\big\}
\label{ erw112f }\\
&&  -\sum_{j=1}^{N}\sum_{k=1}^{2^{j-1}} 
\#\Big\{S\big((2k-2)2^{n-j}, (2k-1)2^{n-j}\big]\cap S\big((2k-1)2^{n-j}, 
(2k)2^{n-j}\big]\Big\}. \nonumber
\end{eqnarray}
Setting 
$$ \beta_k={\#\big\{S((k-1)2^{n-N}, k2^{n-N}]\big\}}  $$ and 
$$\alpha_{j,k}= { \#\Big\{S\big((2k-2)2^{n-j}, (2k-1)2^{n-j}\big]\cap 
S\big((2k-1)2^{n-j}, (2k)2^{n-j}\big]\Big\} }$$ 
 leads to the 
decomposition 
$$ 
\ol{R}_{2^n}=\sum_{k=1}^{2^N}\ol \beta_k-\sum_{j=1}^N 
\sum_{k=1}^{2^{j-1}}\ol \alpha_{j,k}. 
$$ Recall that (Lemma 3 in \cite{Cb}), 
\begin{equation} 
\sup_n\E\exp\Big\{\lambda{\log n\over 
n}\#\big\{S[1,n]\big\}\Big\}<\infty 
\label{erw113} 
\end{equation} for all $\lambda>0$. In particular, 
$$ 
\sup_n\E\exp\Big\{\lambda{\log 2^{n-N}\over 
2^{n-N}}\vert\ol\beta_1\vert\Big\} <\infty. 
$$ Notice that $\ol\beta_1,\cdots,\ol\beta_{2^N}$ is an i.i.d.\ sequence 
   with $\E\ol\beta_1=0$. By Lemma 1 in \cite{BCR}, there is a $\theta >0$ 
such that 
$$ 
\sup_n\E\exp\Big\{\theta 2^{-N/2}{\log 2^{n-N}\over 2^{n-N}}\Big\vert 
\sum_{k=1}^{2^N}\ol\beta_k\Big\vert\Big\}<\infty. 
$$ By the choice of $N$ one can see that there is a $c>0$ independent of 
$n$ such that 
$$ 2^{-N/2}{\log 2^{n-N}\over 2^{n-N}}\ge c{(\log 2^{n})^2\over 2^{n}}. 
$$ So there is some $\theta >0$ such that 
$$ 
\sup_n\E\exp\Big\{\theta {(\log 2^{n})^2\over 2^{n}}\Big\vert 
\sum_{k=1}^{2^N}\ol\beta_k\Big\vert\Big\}<\infty. 
$$ 
 
We need to show that for some $\theta >0$, 
\begin{equation} 
\sup_n\E\exp\bigg\{\theta {(\log 2^n)^2\over 2^{n}}\Big\vert 
\sum_{j=1}^N\sum_{k=1}^{2^{j-1}}\ol\alpha_{j,k} 
\Big\vert\bigg\} <\infty.\label{erw114} 
\end{equation} 
 
Set 
\begin{equation} \tilde J_n=\#\big\{S[1,n]\cap S'[1,n]\big\}\hskip.2in 
n=1,2,\cdots,\label{erw115} 
\end{equation} where $S'$ is an independent copy of the random walk 
$S$. In our notation, for each $1\le j\le N$, 
$\{\ol\alpha_{j,1},\cdots, \ol\alpha_{j,2^{j-1}}\}$ is an i.i.d. sequence 
with the same distribution as $\tilde J_{2^{n-j}}$. 
   By Theorem \ref{rd12} (with $p=2$), there is a $\delta>0$ such that 
$$ 
\sup_n\sup_{j\leq N}\E\exp\Big\{\delta {(\log 2^{n-j})^2\over 
2^{n-j}}\big\vert 
\ol\alpha_{j,1} 
\big\vert\Big\}<\infty. 
$$ 
    By Lemma 1 in 
\cite{BCR} again, there is a $\bar\theta >0$ such that 
$$\sup_{n}\sup_{j\leq N}\E\exp\bigg\{\bar\theta 2^{-j/2}{(\log 
2^{n})^2\over 2^{n-j}} 
\Big\vert 
\sum_{k=1}^{2^{j-1}}\ol\alpha_{j,k} 
\Big\vert\bigg\} <\infty.$$ Hence for some $\theta >0$ 
$$ C(\theta) 
\equiv \sup_{n}\sup_{j\leq N}\E\exp\bigg\{\theta 2^{j/2}{(\log 2^{n})^2 
\over 2^n} 
\Big\vert 
\sum_{k=1}^{2^{j-1}}\ol\alpha_{j,k} 
\Big\vert\bigg\}<\infty. 
$$ 
 
Write 
$$ 
\lambda_N= \prod_{j=1}^N\big(1-2^{-j/2}\big)\hskip.1in\hbox{and} 
\hskip.1in \lambda_\infty = \prod_{j=1}^\infty \big(1-2^{-j/2}\big). 
$$ 
Using H\"older's inequality with $1/p=1-2^{-N/2},\,1/q=2^{-N/2}$ we 
have 
\begin{align*} 
\E\exp\bigg\{&\lambda_N\theta {(\log 2^n)^2\over 2^{n}}\Big\vert 
\sum_{j=1}^N\sum_{k=1}^{2^{j-1}}\ol\alpha_{j,k} 
\Big\vert\bigg\}\\ &\le \bigg(\E\exp\bigg\{\lambda_{N-1}\theta {(\log 
2^n)^2\over 2^{n}}\Big\vert 
\sum_{j=1}^{N-1}\sum_{k=1}^{2^{j-1}}\ol\alpha_{j,k} 
\Big\vert\bigg\}\bigg)^{1-2^{-N/2}}\\ &\hskip.3in\times 
\bigg(\E\exp\bigg\{\lambda_N \theta 2^{N/2}{(\log 2^n)^2 
\over 2^{n}} 
\Big\vert 
\sum_{k=1}^{2^{N-1}}\ol\alpha_{N,k}\Big\vert\bigg\}\bigg)^{2^{-N/2}}\\ 
&\le \E\exp\bigg\{\lambda_{N-1}\theta {(\log 2^n)^2\over 
2^{n}}\Big\vert 
\sum_{j=1}^{N-1}\sum_{k=1}^{2^{j-1}}\ol\alpha_{j,k} 
\Big\vert\bigg\}\cdot C(\theta)^{2^{-N/2}} 
\end{align*}
since $ \lambda_N<1  $. 
Repeating this procedure, 
\begin{align*} 
\E\exp\bigg\{&\lambda_N\theta {(\log 2^n)^2\over 2^n}\Big\vert 
\sum_{j=1}^N\sum_{k=1}^{2^{j-1}}\ol\alpha_{j,k} 
\Big\vert\bigg\}\\ &\le C(\theta)^{2^{-1/2}+\cdots +2^{-N/2}} 
\le C(\theta)^{2^{-1/2}(1-2^{-1/2})^{-1}}. 
\end{align*}  So we have 
$$ 
\sup_n\E\exp\bigg\{\lambda_\infty\theta {(\log 2^n)^2\over 
2^n}\Big\vert 
\sum_{j=1}^N\sum_{k=1}^{2^{j-1}}\ol\alpha_{j,k} 
\Big\vert\bigg\}\le C(\theta)^{2^{-1/2}(1-2^{-1/2})^{-1}}. 
$$ We have proved (\ref{erw114}) and therefore (\ref{erw112}) when 
$n$ is the power of 
$2$. 
\medskip We now prove Theorem \ref{rd13} for general $n$. Given an 
integer 
$n\ge 2$, we have the following unique representation: 
$$ n=2^{m_1}+2^{m_2}+\cdots +2^{m_l} 
$$ where $m_1>m_2>\cdots m_l\ge 0$ are integers. Write 
$$ n_0=0\hskip.1in \hbox{and}\hskip .1in n_i=2^{m_1}+\cdots +2^{m_i} 
\hskip.2in i=1,\cdots, l. 
$$ Then 
\begin{align*} 
\#\big\{S[1,n]\big\}&=\sum_{i=1}^l\#\big\{S(n_{i-1}, n_i]\big\} 
-\sum_{i=1}^{l-1}\#\big\{S(n_{i-1}, n_i]\cap S(n_i,n]\big\}\\ 
&=\sum_{i=1}^lB_i-\sum_{i=1}^{l-1}A_i. 
\end{align*} 
  Write 
$$ 
\sum_{i=1}^lB_i={\sum_i}'B_i+{\sum_i}''B_i 
$$ where ${\sum_i}'$ is the summation over $i$ with 
$2^{m_i}\ge\sqrt{n}$ and ${\sum_i}''$ is the summation over $i$ with 
$2^{m_i}<\sqrt{n}$. We also define the products ${\prod_i}'$ and 
${\prod_i}''$  in a similar  manner. Then 
\begin{align*} 
\E\exp\Big\{\theta& {(\log n)^2\over n}\Big\vert {\sum_i}'(B_i-\E 
B_i)\Big\vert\Big\}\\ &\le {\prod_i}'\bigg(\E\exp\Big\{\theta {(\log 
n)^2\over n}2^{-m_i} 
\Big({\sum_j}'2^{m_j}\Big)\vert\ol{R}_{2^{m_i}}\vert\Big\} 
\bigg)^{2^{m_i}\Big({\sum_j}'2^{m_j}\Big)^{-1}}\\ &\le 
{\prod_i}'\bigg(\E\exp\Big\{4\theta {(\log 2^{m_i})^2\over 2^{m_i}} 
\vert\ol{R}_{2^{m_i}}\vert\Big\} 
\bigg)^{2^{m_i}\Big({\sum_j}'2^{m_j}\Big)^{-1}}\\ &\le 
\sup_m\E\exp\Big\{4\theta {(\log 2^m)^2\over 2^m} 
\vert\ol{R}_{2^m}\vert\Big\}. 
\end{align*} 
 
Assume that the set $\{1\le i\le l;\hskip.05in 2^{m_i}<\sqrt{n}\}$ is 
non-empty. We have 
$$ {\sum_i}'' 2^{m_i}\le 2\sqrt{n}. 
$$ So we have 
$$ {(\log n)^2\over n}\le {1\over\sqrt{n}}\le 
2\Big({\sum_i}''2^{m_i}\Big)^{-1}. 
$$ Hence 
\begin{align*} 
\E\exp\Big\{&\theta {(\log n)^2\over n}\Big\vert {\sum_i}''(B_i-\E 
B_i)\Big\vert\Big\}\\ &\le {\prod_i}''\bigg(\E\exp\Big\{\theta {(\log 
n)^2\over n}2^{-m_i} 
\Big({\sum_j}''2^{m_j}\Big)\vert\ol{R}_{2^{m_i}}\vert\Big\} 
\bigg)^{2^{m_i}\Big({\sum_j}''2^{m_j}\Big)^{-1}}\\ &\le 
{\prod_i}''\bigg(\E\exp\Big\{2\theta {1\over 2^{m_i}} 
\vert\ol{R}_{2^{m_i}}\vert\Big\} 
\bigg)^{2^{m_i}\Big({\sum_j}''2^{m_j}\Big)^{-1}}\\ &\le 
\sup_m\E\exp\Big\{2\theta {1\over 2^m} 
\vert\ol{R}_{2^m}\vert\Big\}. 
\end{align*} By the Cauchy-Schwarz inequality and what we have proved 
in the previous step, there exists $\theta>0$ such that 
$$ 
\E\exp\Big\{\theta {(\log n)^2\over n}\Big\vert 
\sum_{i=1}^l(B_i-\E B_i)\Big\vert\Big\} 
$$ is bounded uniformly in $n$. 
\medskip By the fact that 
\begin{equation} n-n_i=2^{m_{i+1}}+\cdots +2^{m_l}\le 
2^{m_i}\label{1.15s} 
\end{equation} 
  we have 
\begin{equation} 
  A_i\buildrel d\over =\#\big\{S[1,2^{m_i}]\cap S'[1, n-n_i]\big\} 
\le J_{2^{m_i}}.\label{1.15t} 
\end{equation} 
By (\ref{erw14b}) there is a constant $C>0$ independent of 
$n$ such that 
\bea 
&& 
\sum_{i=1}^{l-1}\E A_i\le \sum_{i=1}^{l-1}\E J_{2^{m_i}} 
\le C\sum_{i=1}^{l}{2^{m_i}\over m_i^2}\label{1.15r}\\ 
&&\hspace{ .5in}\le C\sum_{m_i<l/2}{2^{m_i}\over m_i^2} 
+C\sum_{m_i\geq  l/2}^{l}{2^{m_i}\over m_i^2}\nn\\ 
&&\hspace{ .5in}\le C 2^{l/2} 
+C{n\over (\log n)^2}\nn\\ 
&&\hspace{ .5in}\le C{n\over (\log n)^2}.\nn 
\eea 
It remains to show that 
\begin{equation} 
\sup_n\E\exp\Big\{\theta {(\log n)^2\over n} 
\sum_{i=1}^{l-1} A_i\Big\}<\infty.\label{erw115a} 
\end{equation} 
 
Using (\ref{1.15r},) this follows from (\ref{erw19}), (with $p=2$), and the 
same argument used for 
$B_1-\E B_1,\cdots, B_l-\E B_l$. 
\qed 
 
In view of the remark prior to Theorem \ref{rd13}, the next result shows 
that 
$\ol{R}_n$ has a non-symmetric tail behavior. 
 
\begin{theorem} \label{rd14} Under the assumptions of Theorem 
\ref{rd13}, 
\begin{equation} 
\sup_n\E\exp\Big\{\theta{(\log n)^2\over n}\ol{R}_n\Big\} 
<\infty\label{erw116} 
\end{equation} for all $\theta>0$. 
\end{theorem}

\proof  By Theorem \ref{rd13}, (\ref{erw116}) holds for some 
$\theta_0>0$. For $\theta >\theta_0$, take an integer $m\ge 1$ such that 
$m^{-1}\theta <\theta_0$. It is easy to see that it suffices to prove 
\begin{equation} 
\sup_n\E\exp\Big\{\theta {(\log n)^2\over mn}\ol{R}_{nm}\Big\}<\infty. 
\label{erw117} 
\end{equation} Set $\zeta_{jn}=\#\{S((j-1)n,jn]\}$. 
  By the facts that 
$$ 
\ol{R}_{nm}\le\sum_{j=1}^m \overline{\zeta_{jn}} 
+\Big(\sum_{j=1}^m\E\zeta_{jn}\Big)-\E R_{nm} 
$$ and that by (\ref{secr11}) and (\ref{secr12}), 
\begin{align*} 
\Big(\sum_{j=1}^m\E\zeta_{jn}\Big)-\E R_{mn} &=m\E R_n-\E 
R_{mn}\\
&=\frac{mn}{\HH(n)}+O\Big(\frac{mn}{\HH(n)^2}\Big) -\frac{mn}{\HH(mn)}
+O\Big(\frac{mn}{\HH(mn)^2}\Big)\\
&=\frac{mn}{\HH(n)\HH(mn)}((\HH(mn)-\HH(n)) +O\Big(\frac{n}{\log^2 n}\Big)\\
&=O\Big({n\over(\log n)^2}\Big) 
\end{align*}
 as $n\to \infty$ (note $m$ is fixed), there is a constant $C_{m,\theta}>0$ depending only 
on $m$ and $\theta$ such that 
$$ 
\E\exp\Big\{\theta {(\log n)^2\over mn}\ol{R}_{nm}\Big\} 
\le C_m\bigg(\E\exp\Big\{\theta {(\log n)^2\over 
mn}\ol{R}_{n}\Big\}\bigg)^m. 
$$ So we have (\ref{erw117}). 
\qed

\def\be{{\bf e}} 
 
\section{Moderate deviations for $R_n-\E R_n$} 
 
We can now prove Theorem \ref{theo-range}. 
\medskip

\proof 
We first prove the upper bound. Let $t>0$ and write 
$K=[t^{-1}b_n]$. Divide $[1,n]$ into $K$ disjoint subintervals, each of 
length $[n/K]$ or $[n/K]+1$. Call the $i^{th}$ subinterval $I_i$. Let 
$E_i=\#\{S(I_i)\}$.  Then 
$$ 
\ol{R}_n\le \sum_{j=1}^K\ol{E}_j+\Big(\sum_{j=1}^K\E {E}_j\Big)-\E 
R_n 
$$ 
From (\ref{lrth6.9}) we have 
\begin{align}
\sum_{j=1}^K&\E E_j-\E R_n\label{secr2}\\ 
&=K{n/K\over \HH([n/K])}-{n\over \HH(n)}
+{1\over 2\pi\sqrt{\det\Gamma}}\Big\{K{n/K\over \HH^2([n/K])}
-{n\over \HH^2(n)}\Big\}
+o\Big({n\over \HH^2(n)}\Big)\nn\\
&={n(\HH(n)-\HH([n/K]))\over \HH^2(n)}\Big\{1+{\HH(n)-\HH([n/K])\over \HH([n/K])}\Big\}\nn\\
&~~~~+{n\over \HH^2(n)}\Big\{
{\HH^2(n)-\HH^2([n/K])\over \HH^2([n/K])}\Big\}+o\Big({n\over \HH^2(n)}\Big),\nn
\end{align} 
where the error term can be taken to be independent of 
$\{b_n\}$. 
(This is where the hypothesis $\log b_n=o((\log n)^{1/2})$ is used.)
Since
$$
\HH(n)-\HH([n/K])=\sum_{k=[n/K]+1}^n \P\{S_k=0\}
\sim {\log K\over 2\pi\sqrt{\det\Gamma}},
$$
we have
\begin{equation}\label{secr3}
\sum_{j=1}^K\E E_j-\E R_n
={n(\HH(n)-\HH([n/K]))\over \HH^2(n)}
+o\Big({n\over \HH^2(n)}\Big).
\end{equation}

Hence for any $\lambda>0$,  
\begin{align*}
\P\Big\{\ol{R}_n&\ge         
{n\over \HH^2(n)}\big(\HH(n)-\HH([n/b_n])\big)\Big\}\\
&\le
\exp\Big\{-\lambda b_n\big(\HH(n)-\HH([n/b_n])\big)\Big\}
\E\exp\Big\{\lambda {\HH^2(n)b_n\over n}\ol{R}_n\Big\}\\
&\le\exp\Big\{-\lambda b_n\big(\HH([n/K])-\HH([n/b_n])\big)+o(b_n)\Big\}
\bigg(\E\exp\Big\{\lambda {\HH^2(n)b_n\over n}\ol{E}_1\Big\} \bigg)^K.
\end{align*}

Notice that
$$
\lim_{n\to\infty}\big(\HH([n/K])-\HH([n/b_n])\big)
={\log t\over 2\pi\sqrt{\det\Gamma}}
$$
and that
by \cite[Theorem 6.1]{LeGall6}, 
$$
{\HH^2(n)b_n\over n}\ol{E}_1\buildrel d\over\longrightarrow
-{2\pi t\over 2\pi\sqrt{\det\Gamma}}\gamma_1,
$$
where $\gamma_t$ is the renormalized self-intersection local time of a 
planar Brownian motion. By Theorem 
\ref{rd14} and the dominated convergence theorem, 
$$ 
\E\exp\Big\{\lambda {\HH^2(n)b_n\over n}\ol{E}_1\Big\}
\longrightarrow\E\exp\Big\{-\lambda{2\pi t\over 2\pi\sqrt{\det\Gamma}}\gamma_1
\Big\}
$$

Consequently,
\begin{align}
\limsup_{n\to\infty}&\, b_n^{-1}\log\P\Big\{\ol{R}_n\ge         
{n\over \HH^2(n)}\big(\HH(n)-\HH([n/b_n])\big)\Big\}\label{secr1}\\
&\le -\lambda{\log t\over 2\pi\sqrt{\det\Gamma}}
+{1\over t}\log \E\exp\Big\{-\lambda{2\pi t\over 2\pi\sqrt{\det\Gamma}}\gamma_1
\Big\}\nn\\
&={\lambda\over2\pi\sqrt{\det\Gamma}}\log{\lambda\over2\pi\sqrt{\det\Gamma}}\nn\\
&~~~+{1\over t}\log \E\exp\Big\{-{\lambda t\over2\pi\sqrt{\det\Gamma}}\log
{\lambda t\over2\pi\sqrt{\det\Gamma}}
-\lambda{2\pi t\over 2\pi\sqrt{\det\Gamma}}\gamma_1
\Big\}.\nn   
\end{align}
By \cite{BC}, see the proof of Theorem 3.2, the limit
\begin{equation}
C\equiv\lim_{t\to\infty}{1\over t}\log\E\exp\big\{- t\log t-2\pi 
t\gamma_1\big\}\label{pat.1 }
\end{equation} exists. 
Set
$$L=\exp(-1-C).$$
Letting $t\to\infty$ in (\ref{secr1}) gives
\begin{align*}
\limsup_{n\to\infty}&b_n^{-1}\log\P\Big\{\ol{R}_n\ge         
{n\over \HH^2(n)}\big(\HH(n)-\HH([n/b_n])\big)\Big\}\\
&\le {\lambda\over2\pi\sqrt{\det\Gamma}}\log{\lambda\over2\pi\sqrt{\det\Gamma}}
+C{\lambda\over2\pi\sqrt{\det\Gamma}}.
\end{align*}
Taking
$$
{\lambda\over2\pi\sqrt{\det\Gamma}}=\exp\big\{-1-C\big\}
$$
then yields
\begin{align*}
\limsup_{n\to\infty}&b_n^{-1}\log\P\Big\{\ol{R}_n\ge         
{n\over \HH^2(n)}\big(\HH(n)-\HH([n/b_n])\big)\Big\}\\
&\le -\exp\big\{-1-C\big\}=-L.
\end{align*}

\bigskip

We now prove the lower bound. 
    The proof is similar to that of Proposition 4.4 of \cite{BK}. 
Fix $n$ and let 
$K=[b_n]$. Let $M=[n/b_n]$. Let $I_j$ be the interval 
$(m_j,m_{j+1}]$, where the  $m_j$ are integers such that $m_0=0$, 
$m_K=n$, and $m_{j+1}-m_j$ is equal to either 
$M$ or $M+1$. 
 
Let $\be$ be a vector of length $\sqrt M$ and let $B(x,r)$ be the ball of 
radius $r$ about $x$. Set 
$$E_j=\# \{S(I_j)\}, \qquad H_j=\#\{S(I_j)\cap S(I_{j-1})\}.$$  Let  
\begin{equation} A_j=\{S_{m_{j+1}}\in B((j+1)\be, \tfrac18 \sqrt M)\}\cap 
\{ S(I_j)\subset B((j+\tfrac12)\be, \sqrt M)\}\end{equation} 
and 
\begin{equation} 
B_j=\{ \ol E_j \log^2 M/M\geq -c_1\}\label{erd25a} 
\end{equation} 
where we will select $c_1$ in a moment. By the central limit 
theorem, we know 
$\P^{S_{m_{j-1}}}(A_j)\geq c_2$ on the event $A_{j-1}$ if $n$ is large. 
By \cite[Theorem 6.1]{LeGall6}, 
$\P^{S_{m_{j-1}}}(A_j\cap B_j)>c_2/2$  on the event $A_{j-1}$ if we take 
$c_1$ sufficiently large. 
    If we let 
$$F=\bigcap_{j=0}^{K-1} (A_j\cap B_j),$$ then by the Markov property 
applied 
$K-1$ times we have 
\begin{equation} \P(F)\geq (c_2/2)^{K-1}.\label{erd26} \end{equation} 
 
    On the set $F$ we have that 
$S(I_j)$ is disjoint from $S(I_i)$ if $|i-j|>1$, and so on $   F$ 
\begin{equation} \ol R_n=\sum_{j=1}^K \ol E_j+\Big( 
\Big(\sum_{j=1}^K \E E_j\Big)-\E R_n\Big)- 
\sum_{j=1}^K H_j.\label{erd27} \end{equation} On the set $F$ the event 
$B_j$ holds for each $j$, and so 
\begin{equation} \sum_{j=1}^K \ol E_j\geq -\frac{c_1KM}{\log^2 
M}\geq -\frac{c_3 n}{\log^2 n}. 
\label{erd28} \end{equation} As in (\ref{secr3}),
\begin{equation}\label{erw39}
 \Big(\sum_{j=1}^K \E E_j\Big)-\E R_n  
=\frac{n(\HH(n)-\HH([n/K]))}{\HH(n)^2}+o\Big(\frac{n}{\HH(n)^2}\Big)
\end{equation} if $n$ is large. 
 
\medskip Let $\Lambda >0$ be chosen in a moment. Let  
$$  C_1=\lc \sum_{\{j \ odd\}} H_j\geq \frac{n\Lambda }{\log^2 n}\rc,  
\hspace{.3 in}C_2=\lc \sum_{\{j \ even\}} H_j\geq \frac{n\Lambda }{\log^2 n}\rc.$$  
Set 
$G=F\cap C_1^c\cap C_2^c$. For $j$ odd the $H_j$ are independent, 
and by Lemma 4.6 of \cite{BK}, 
\begin{align*} \P(C_1) &=\P\Big( \sum_{\{j \ odd\}} \frac{H_j}{M/\log^2 
M}\geq c_4K\Lambda \Big)\\ &\leq e^{-c_4c_5K\Lambda } \E e^{c_5 H_j \log^2 M/M}\\ 
&\leq e^{-c_4c_5 K\Lambda } c_6^K, 
\end{align*} where $c_4, c_5, c_6$ do not depend on $\Lambda $ and without 
loss of generality we may assume $c_6>1$. Choose $\Lambda $ large so that 
$e^{-c_4c_5\Lambda }\leq c_6^{-2}$. When $n$ is large, $K$ will be large, and 
then $\P(C_1)\leq \P(F)/3$. We have a similar estimate for $\P(C_2)$, so 
$$\P(G)\geq (c_2/2)^{K-1}/3.$$ 
Set $v_n=\HH(n)-\HH([n/b_n])$.
On the event $G$ 
\begin{equation} \sum_{j=1}^K H_j \leq 2\frac{n\Lambda }{\log^2 
n},\label{erd210}\end{equation} and so combining (\ref{erd28}), 
(\ref{erw39}), and (\ref{erd210}), on the event $G$ 
\begin{equation} \ol R_n\geq \Big(1-\frac{c_7}{v_n}\Big) n 
v_n/\HH( n)^2.\label{erd211}\end{equation} Therefore 
\begin{equation} 
\P\Big(\ol R_n\geq \Big(1-\frac{c_7}{v_n}\Big) n
v_n/\HH( n)^2\Big)\geq c_8 c_9^{b_n}.\label{rwe212} 
\end{equation} Define $b'_n$ by $v'_n=\HH(n)-\HH([n/b'_n])= v_n+c_7$. If we apply 
(\ref{rwe212}) with 
$b_n$ replaced by $b'_n$, we have 
\begin{align*} 
\P(\ol R_n\geq & nv_n/\HH( n)^2)\\ &=\P\Big(\ol  R_n\geq 
\Big(1-\frac{c_7}{ v_n'}\Big) nv'_n/\HH( n)^2\Big)\\ 
&\geq c_8 c_9^{b'_n}. 
\end{align*} We now take the logarithms of both sides, divide by $b_n$, 
and use the fact that the ratio $b_n/b_n'$ is
bounded above and below by positive constants  to 
obtain the lower bound. 
\qed

\noindent {\bf Proof of Corollary \ref{corrange}:}
Assume first that $S_n$ is strongly aperiodic.
We have by Proposition \ref{p21a} that
\begin{equation}\label{1.2}
\P(S_n=0)={ 1\over 2\pi n\sqrt{\det
\Gamma}} +O\Big({1 \over n(\log n)^{1/2}}\Big).
\end{equation}
Then, if $\ga$ denotes Euler's constant
\begin{equation}
\sum_{k=1}^n { 1\over k}=\log n+\gamma+O\Big({ 1\over
n}\Big)\label{1.2e}
\end{equation}
and
\begin{equation}\label{1.2e2}
\sum_{k=3}^n{1 \over k(\log
k)^{1/2}}\leq \int_2^n \frac{dx}{x(\log x)^{1/2}}\leq 
 c_1(\log n)^{1/2}
\end{equation}
so that
\begin{align}
  \HH(n)&=\sum_{k=0}^n \P^0(S_k=0)=1+\frac{1}{2\pi
\sqrt{\det \Gamma}}\sum_{k=1}^n\({ 1\over k}+O\Big({1 \over k(\log
k)^{1/2}}\Big)\) \label{1.3}\\
&=\frac{1}{2\pi
\sqrt{\det \Gamma}}\(\log n+\gamma+O\Big( (\log
n)^{1/2}\Big)\)\nn\\
&=\frac{\log n}{2\pi
\sqrt{\det \Gamma}}\(1+O\Big({1 \over (\log
n)^{1/2}}\Big)\).\nn
\end{align}
Similarly we 
\begin{equation}\label{1.3k}
\sum_{k=[n/b_n]+1}^{ n}{1 \over
k(\log k)^{1/2}}\leq 
c_2\((\log n)^{1/2}-(\log (n/b_n))^{1/2}\).
\end{equation}
To evaluate this note that
\begin{align} 
(\log (n/b_n))^{1/2}&=(\log n-\log b_n)^{1/2}\label{1.3l}\\ &
=(\log n)^{1/2}(1-\log b_n/\log n)^{1/2}\nonumber\\ &
=(\log n)^{1/2}(1+ O(\log b_n/ \log n))\nonumber\\ &
=(\log n)^{1/2}+O(\log b_n/(\log n)^{1/2})
\nonumber
\end{align}
by our assumption that $\log  b_n=o((\log n)^{1/2})$.
It follows that
\begin{align}
  \HH(n)-\HH([n/b_n])&=\frac{1}{2\pi
\sqrt{\det \Gamma}}\sum_{k=[n/b_n]+1}^n\({ 1\over k}+O\Big({1 \over
k(\log  k)^{1/2}}\Big)\) \label{1.4}\\
&=\frac{1}{2\pi
\sqrt{\det \Gamma}}\(\log b_n+O\Big({\log b_n \over (\log
n)^{1/2}}\Big)\)\nn\\
&=\frac{\log b_n}{2\pi
\sqrt{\det \Gamma}}\(1+O\Big({1 \over (\log
n)^{1/2}}\Big)\).\nn
\end{align}
We then have that
\begin{align} 
{n\over\HH( n)^2}&(\HH(n)-\HH([n/b_n]))\label{1.5}\\ &
=2\pi
\sqrt{\det \Gamma}\frac{n\log b_n}{(\log n)^{ 2}}\(1+O\Big({1 \over (\log
  n)^{1/2}}\Big)\)\nonumber\\ &
=2\pi
\sqrt{\det \Gamma}\frac{n\log b_n}{(\log n)^{ 2}}(1+a_{ n}),\nonumber
\end{align}
where we use the last equality to define $a_n$.
Let
\begin{equation}
1+\wh a_{ n}=(1+a_{ n})^{ -1}=1+O\Big({1 \over (\log
  n)^{1/2}}\Big).\label{1.5a}
\end{equation}
Then if we set
\begin{equation}
\wh b_{ n}=: b_n^{1+\wh a_{ n} }= b_n^{(1+a_{ n})^{ -1}}\label{}
\end{equation}
we see from (\ref{1.5}) that
\begin{equation}
{n\over\HH( n)^2}(\HH(n)-\HH([n/\,\,\wh b_n]))=
2\pi
\sqrt{\det \Gamma}\frac{n\log b_n}{(\log n)^{ 2}}.\label{1.5b}
\end{equation}
Also, $\log \wh b_n= ( 1+\wh a_{ n})\log b_n=o((\log n)^{1/2})$, so that
Theorem
\ref{theo-range} applies to
  $\wh b_n$, and indeed to $\wh b_n^{ \th}$ for any $\th>0$.

Note that
\begin{eqnarray}
  \wh b_n^{ \th}&=&b_n^{\th\(1+O\Big({1 \over (\log
n)^{1/2}}\Big)\) }\label{1.6}\\ && = b_n^{ \th} \exp\(O\Big({ \log
b_n\over (\log  n)^{1/2}}\Big)\)\nonumber\\ && = b_n^{ \th} (1+o(
1_{ n}))\nonumber
\end{eqnarray}
by our assumption that $\log  b_n=o((\log n)^{1/2})$.    Hence by
(\ref{1.5b}) and (\ref{1.6})
\begin{align} 
\wh b_n^{-\th}\log
&\P\Big\{\ol{R}_n
\ge {\th n\over\HH( n)^2}(\HH(n)-\HH([n/\wh b_n]))\Big\}\label{1.6a}\\ &
= (1+o( 1_{ n}))b_n^{-\th}\log
\P\Big\{\ol{R}_n
\ge 2\pi \th
\sqrt{\det \Gamma}\frac{n\log b_n}{(\log n)^{ 2}}\Big\}\nonumber
\end{align}

Together with Proposition \ref{p21a}, Theorem \ref{theo-range} applied to
  $\wh b_n^{ \th}$ proves the corollary in the strongly aperiodic case.
The modifications to handle the case where $S_n$ is not strongly aperiodic
are very similar to those in Section 2 of \cite{LR}.
\qed

\section{Moderate deviations for $\E R_n -R_n$} 
 
To avoid difficulties connected with subdividing time intervals, it is more 
convenient to look at the continuous time analogue of 
$S_n$. We let $T_1, T_2, \ldots$ be i.i.d.~ exponential random variables 
with parameter 1 that are independent of the  sequence $S_n$.  Define 
$Z_t=S_n$ if 
$\sum_{i=1}^n T_i\leq t< \sum_{i=1}^{n+1} T_i$. 
$Z_t$ is { a L\'evy} process that waits an exponential  length of 
time, then jumps according to $X_1$, and then repeats the procedure. 
Define $N_t=n$ if 
$\sum_{i=1}^n T_i\leq t< \sum_{i=1}^{n+1} T_i$.  Note that $N_t$ is a 
Poisson process with $\E N_t=t$ and that $Z_t=S_{N_t}$.  
We write $|Z[a,b]|$ for the cardinality of 
$\{Z_s: s\in [a,b]\}$. 
 
Theorems \ref{rd12} and \ref{rd13} have the following analogues for 
continuous time processes. We omit the proofs, which are almost identical to 
the proofs given for the discrete time random walks. 
 
\begin{lemma}\label{LU2} Let $Z_1(t),\cdots, Z_p(t)$ be independent 
copies of 
$Z(t)$. There is $C>0$ such that 
\begin{equation} 
\label{4.23} 
\sup_{y_1,\cdots,y_p}\E^{(y_1,\cdots, y_p)}\Big\vert 
Z_1[0,t]\cap\cdots\cap Z_p[0,t]\Big\vert^m 
\le C^m (m!)^{p-1}\Big({t\over (\log t)^p}\Big)^m. 
\end{equation} 
Consequently, there is $\theta >0$ such that 
\begin{equation} 
\label{4.24} 
\sup_t\sup_{y_1,\cdots,y_p}\E^{(y_1,\cdots, y_p)}\exp\bigg\{\theta\Big({(\log t)^p\over t} 
\Big\vert Z_1[0,t]\cap\cdots\cap 
Z_p[0,t]\Big\vert\Big)^{(p-1)^{-1}}\bigg\} <\infty. 
\end{equation} 
\end{lemma}

\begin{lemma}\label{LU3} There is $\theta >0$ such that 
\begin{equation} 
\label{4.25} 
\sup_t\E\exp\Big\{\theta{(\log t)^2\over t}\big\vert\E\vert Z[0,t]\vert 
-\vert Z[0,t]\vert\big\vert\Big\}<\infty. 
\end{equation} 
Consequently 
\begin{equation} 
\limsup_{t\to\infty}{1\over b_t}\log\P\Big\{\Big\vert 
\E\big\vert Z[0,t]\big\vert -\big\vert Z[0,t]\big\vert\Big\vert\ge 
\lambda {tb_t\over (\log t)^2}\Big\} 
\le -\theta\lambda. 
\label{4.123} 
\end{equation} 
\end{lemma}

We will prove Theorem \ref{LU1.3} by first proving the following analogue 
for 
$Z_t$. 
 
\begin{theorem}\label{TLU1} For any $\lambda >0$ and for any $b_t$ 
satisfying 
$ b_t\to \infty$ and 
$b_t=o\Big((\log t)^{1/5}\Big)$ as $t\to\infty$, we have 
\begin{eqnarray} && 
\label{4.1} 
\lim_{t\to\infty}{1\over b_t}\log\P\Bigg\{\bigg\vert\E\big\vert 
Z[0,t]\big\vert -\big\vert Z[0,t]\big\vert\bigg\vert\ge \lambda {tb_t\over 
(\log t)^2}\Bigg\} 
\\ && =-(2\pi)^{-2}\det(\Gamma)^{-1/2}\kappa 
(2,2)^{-4}\lambda.\nonumber 
\end{eqnarray} 
\end{theorem} 
 
\medskip The next proposition shows that Theorem \ref{LU1.3} follows from 
Theorem 
\ref{TLU1} and Theorem \ref{theo-range}. 
 
\begin{proposition}\label{LU0} For any $\eps>0$, 
\begin{equation} 
\label{4.2} 
\lim_{n\to\infty}{1\over b_n}\log\P\Big\{\Big\vert 
\ov {\vert Z[0,n]\vert}- \ov {\vert S[0,n]\vert} 
\Big\vert\ge \eps {nb_n\over 
(\log n)^2}\Big\} =-\infty. 
\end{equation} 
\end{proposition} 
 
\begin{remark}{\rm Our proof actually gives a stronger result, but this is all 
we need.} 
\end{remark} 
 
\proof Observe that if $n>m$, then 
\begin{equation} 
 \Big\vert\E|S[0,n]|-\E|S[0,m]|\Big\vert\, \leq \E\vert S[m, n]\vert =\E\vert S[0, 
n-m]\vert\le n-m.\label{4.3} 
\end{equation} 
  Consequently, 
\begin{eqnarray} 
\Big\vert\E \vert Z[0, n]\vert -\E\vert S[0, n]\vert\Big\vert  &=&\Big\vert\E 
|S[0, N_n]| -\E|S[0,n]|\,\Big\vert\label{4.4}\\ 
&\le& \E\vert N_n-n\vert\le C\sqrt{n}.  \nonumber 
\end{eqnarray} 
  Hence, it suffices to show that for any $\eps >0$ 
\begin{equation} 
\lim_{n\to\infty}{1\over b_n}\log\P\Big\{\Big\vert 
\vert Z[0,n]\vert -\vert S[0,n]\vert 
\Big\vert\ge \eps{\sqrt{n}b_n^{3/2}\over\log n}\Big\} =-\infty.\label{4.4A} 
\end{equation}

Let $M>0$ be fixed. On the event $\{\vert N_n-n\vert\le M\sqrt{nb_n}\}$ 
\begin{align} 
\Big\vert\vert Z[0, n]&\vert -\vert S[0, n]\vert\Big\vert 
\le \big\vert S\big[N_n\wedge n, N_n\vee n\big]\big\vert\label{4.5}\\ 
&\buildrel d\over =\big\vert S\big[0,\hskip.05in N_n\vee n -N_n\wedge 
n\big]\big\vert 
\le \big\vert S\big[0,\hskip.05in 2M\sqrt{n b_n}\big]\big\vert.\nn 
\end{align} So we have 
\begin{eqnarray} 
&&\P\Big\{\Big\vert 
\vert Z[0,n]\vert -\vert S[0,n]\vert 
\Big\vert\ge \eps{\sqrt{n}b_n^{3/2}\over\log n}\Big\}\nn\\ &&~~~\le 
\P\Big\{\big\vert S\big[0,\hskip.05in 2M\sqrt{n b_n}\big]\big\vert 
\ge\eps{\sqrt{n}b_n^{3/2}\over\log n}\Big\} +\P\big\{\vert 
N_n-n\vert\ge M\sqrt{nb_n} \big\}.\label{4.6} 
\end{eqnarray} 
 
By Lemma 3 in \cite{Cb}, 
\begin{equation} 
\label{4.8} 
\sup_n\E\exp\Big\{\theta {\log n\over \sqrt{nb_n}} S\big[0,\hskip.05in 
2M\sqrt{n b_n}\big]\Big\}<\infty, \qquad \theta >0. 
\end{equation} 
  By the Chebyshev inequality one can see that 
\begin{equation} 
\label{4.9} 
\lim_{n\to\infty}{1\over b_n}\log 
\P\Big\{\big\vert S\big[0,\hskip.05in 2M\sqrt{n b_n}\big]\big\vert 
\ge\eps{\sqrt{n}b_n^{3/2}\over\log n}\Big\}=-\infty. 
\end{equation} 
By the classical moderate deviation principle (\cite[Theorem 3.7.1]{DZ}), 
\begin{equation} 
\label{4.10} 
\lim_{n\to\infty}{1\over b_n}\log 
\P\big\{\vert N_n-n\vert\ge M\sqrt{nb_n} \big\}=-{M^2\over 2}. 
\end{equation} 
 
Thus, 
\begin{equation} 
\label{4.11} 
\limsup_{n\to\infty}{1\over b_n}\log\P\Big\{\Big\vert 
\vert Z[0,n]\vert -\vert S[0,n]\vert 
\Big\vert\ge \eps{\sqrt{n}b_n^{3/2}\over\log n}\Big\} 
\le -{M^2\over 2}. 
\end{equation} 
Letting $M\to\infty$ proves the proposition.\qed 
 
Thus we we need to prove  Theorem 
\ref{TLU1}.  By the G\"artner-Ellis theorem (\cite[Theorem 2.3.6]{DZ}), to prove Theorem 
\ref{TLU1} it suffices to prove 
\begin{eqnarray}&& 
\lim_{t\to\infty}{1\over b_t}\log \E\exp\bigg\{ 
\theta\sqrt{b_t\over t}(\log t) 
\Big\vert \E\big\vert Z[0,t]\big\vert -\big\vert Z[0,t] 
\big\vert\Big\vert^{1/2}\bigg\}\label{4.101}\\ 
&&=(\theta\pi)^2 
\sqrt{\det(\Gamma)}\kappa(2,2)^4.\nn 
\end{eqnarray}

Let $h(x)$ be a smooth symmetric 
probability density on 
$\R^2$ with compact support and write $h_\eps 
(x)=\eps^{-2}h(\eps^{-1} x)$.  We have 
\begin{equation} 
\Lambda_\eps(t)\equiv\sum_{x\in\Z^2}h_\eps 
\Big({x\over\sqrt{t}}\Big) 
\sim t, \qquad t\to \infty. 
\label{4.133n} 
\end{equation} 
 
The following lemma describing  exponential asymptotics for the 
smoothed range will be proved in Section \ref{sec-smoothrange}.

\begin{lemma}\label{LU6} Let 
\begin{equation} 
A_{ t}( \eps)\equiv \Lambda_\eps\Big({t\over 
b_t}\Big)^{-2}\sum_{x\in \Z^2}\bigg[\sum_{y\in Z[0,t]} 
h_\eps\Big(\sqrt{b_t\over t}(x-y)\Big)\bigg]^2.\label{1.00} 
\end{equation} 
For any $\theta >0$, 
\begin{eqnarray} &&\lim_{t\to\infty}{1\over b_t}\log 
\E\exp\bigg\{\theta\sqrt{b_t\over t} 
(\log t)\, 
\vert A_{ t}( \eps) \vert^{1/2}\bigg\}\label{1.00a}\\ &&~~~=\sup_{g\in {\cal 
F}}\bigg\{2\pi\theta\sqrt{\det(\Gamma)} 
\bigg(\int_{\R^2}\vert (g^2\ast h_\eps)(x)\vert^2dx\bigg)^{1/2}\nn\\ && 
\hspace{.3in}-{1\over 2}\int_{\R^2}\langle\nabla g(x),\Gamma\nabla 
g(x)\vert^2 dx\bigg\}.\nn 
\end{eqnarray} 
where 
$$ 
{\cal F}=\{g\in W^{1,2}(\R^2);\hskip.1in \vert\vert g\vert\vert_2 =1\}. 
$$ 
Furthermore, for any $N=0,1,\ldots$ and any $\epsilon >0$, 
\begin{eqnarray} 
&&\lim_{t\to\infty}{1\over b_t}\log 
\E\exp\bigg\{\theta\sqrt{b_t\over t} 
(\log t)\nn\\ 
&&\hskip1in\times\bigg(\Lambda_\eps \Big({t\over 
b_t}\Big)^{-2}\sum_{x\in \Z^2}\bigg[\sum_{y\in Z[0,2^{-N}t]} 
h_\eps\Big(\sqrt{b_t\over t}(x-y)\Big)\bigg]^2 
\bigg)^{1/2}\bigg\}\nn\\ 
&&\hspace{.4in}\leq 2^{-N+2}\pi^2\theta^2\sqrt{\det(\Gamma)} 
\kappa (2,2)^4.\label{4.145} 
\end{eqnarray} 
\end{lemma}

The following lemma on exponential approximation will  be 
proved in Section \ref{sec-expapp}. In this lemma $Z'$ denotes an 
independent copy of $Z$. 
 
\begin{lemma}\label{LU7} 
Let 
\begin{eqnarray}&& 
B^{ ( j)}_{ t}( \eps)\equiv \Lambda_\eps\Big({t\over 
b_t}\Big)^{-2}\nn\\ 
&&\hspace{ .5in}\times \sum_{x\in \Z^2}\bigg[\sum_{y\in Z[0,2^{ -j}t]} 
h_\eps\Big(\sqrt{b_t\over t}(x-y)\Big)\bigg] 
\bigg[\sum_{y'\in Z'[0,2^{ -j}t]} h_\eps\Big(\sqrt{b_t\over 
t}(x-y')\Big)\bigg].\label{1.10} 
\end{eqnarray} 
  Then for any $\theta >0$ and any $j=0,1,\ldots$, 
\begin{eqnarray} 
&& 
\limsup_{\eps\to 0}\limsup_{t\to\infty}{1\over b_t}\nn\\ 
&&\hspace{ .5in}\log 
\E\exp\Big\{\theta\sqrt{b_t\over t}(\log t) 
\big\vert \,\vert Z[0, 2^{ -j}t]\cap Z'[0, 
2^{ -j}t]\vert-B^{ ( j)}_{ t}( \eps)\big\vert^{1/2}\Big\}=0. 
\label{4.81} 
\end{eqnarray} 
\end{lemma} 
 
These lemmas will be the key to proving Theorem \ref{TLU1}. Before 
proving this theorem, we present a simple lemma which will be used several 
times in the proof of Theorem \ref{TLU1}.

\begin{lemma}\label{LU1} Let $l\ge 2$ be a fixed integer and let 
$\{\xi_1(\rho);\hskip.05in \rho>0\},\cdots, 
\{\xi_l(\rho);\hskip.05in \rho>0\}$ be 
$l$ independent non-negative stochastic processes. 
\medskip 
 
(a) If there is a constant $C_1>0$ such that for any $1\le j\le 
l$, 
\begin{equation} 
\label{4.12} 
\limsup_{\rho\to 0^+}\rho\log\P\big\{\xi_j(\rho)\ge\lambda\big\} 
\le -C_1\lambda, \qquad \lambda >0, 
\end{equation} 
  then 
\begin{equation} 
\label{4.13} 
\limsup_{\rho\to 0^+}\rho\log\P 
\Big\{\xi_1(\rho)+\cdots +\xi_l(\rho) 
\ge\lambda\Big\} 
\le -C_1\lambda, \qquad \lambda >0. 
\end{equation} 
\smallskip 
 
(b) If there is a constant $C_2>0$ such that for any $1\le j\le 
l$, 
\begin{equation} 
\label{4.14} 
\limsup_{\rho\to 0^+}\rho\log\E\exp 
\Big\{\rho^{-1}\theta\sqrt{\xi_j(\rho)}\Big\} 
\le C_2\theta^2, \qquad \theta >0, 
\end{equation} then 
\begin{equation} 
\label{4.15} 
\limsup_{\rho\to 0^+}\rho\log\E\exp 
\Big\{\rho^{-1}\theta\sqrt{\xi_1(\rho)+\cdots +\xi_l(\rho)}\Big\} 
\le C_2\theta^2, \qquad \theta >0. 
\end{equation} 
\end{lemma} 
 
\proof. Clearly, part (a) needs only to be proved in the case $l=2$. Given 
$0<\de <\lambda$, let $0=a_0<a_1<\cdots <a_N=\lambda$ be a 
 partition of $[0,\lambda]$ such that $a_k-a_{k-1}<\de$. Then 
\begin{eqnarray} 
\P\big\{\xi_1(\rho)+\xi_2(\rho)\ge \lambda\big\} 
&\le&\sum_{k=1}^N\P\big\{\xi_1(\rho)\in [a_{k-1},a_k]\big\} 
\P\big\{\xi_2(\rho)\ge \lambda -a_k\big\} 
\nn\\ 
&\le&\sum_{k=1}^N\P\big\{\xi_1(\rho)\ge a_{k-1}\big\} 
\P\big\{\xi_2(\rho)\ge \lambda - a_k\big\}. 
\label{4.16} 
\end{eqnarray} 
  Hence 
\begin{eqnarray} 
&& 
\limsup_{\rho\to 0^+}\rho\log 
\P\big\{\xi_1(\rho)+\xi_2(\rho)\ge \lambda\big\}\label{4.17}\\ 
&& 
\le \max_{1\le k\le N}\Big\{-C_1a_{k-1}-C_1(\lambda-a_k)\Big\} 
\le -C_1(\lambda-\de).\nn 
\end{eqnarray} Letting $\de\to 0^+$  proves part (a).

\medskip 
 
We now prove part (b). By Chebyshev's inequality, for any $\lambda >0$ 
\begin{equation} 
\label{4.18} 
\limsup_{\rho\to 0^+}\rho\log 
\P\big\{\xi_j(\rho)\ge \lambda\big\}\le -\sup_{\theta >0} 
\{\theta\sqrt{\lambda} -C_2\theta^2\}=-{\lambda\over 4 C_2}. 
\end{equation} 
  By part (a) 
\begin{equation} 
\label{4.19} 
\limsup_{\rho\to 0^+}\rho\log 
\P\Big\{\xi_1(\rho)+\cdots +\xi_l(\rho)\ge \lambda\big\} 
\le -{\lambda\over 4 C_2}, \qquad \lambda >0. 
\end{equation} 
  In addition, by the triangle inequality and by independence, 
\begin{equation} 
\label{4.20} 
\E\exp 
\Big\{\rho^{-1}\theta\sqrt{\xi_1(\rho)+\cdots +\xi_l(\rho)}\Big\} 
\le\prod_{j=1}^l\E\exp 
\Big\{\rho^{-1}\theta\sqrt{\xi_j(\rho)}\Big\}. 
\end{equation} 
  So by assumption, for any $\theta >0$, 
\begin{equation} 
\label{4.21} 
\limsup_{\rho\to 0^+}\rho\log\E\exp 
\Big\{\rho^{-1}\theta\sqrt{\xi_1(\rho)+\cdots +\xi_l(\rho)}\Big\} 
<\infty. 
\end{equation} 
By   \cite[Lemma 4.3.6]{DZ}, 
\begin{eqnarray} 
\label{4.22} 
&&\limsup_{\rho\to 0^+}\rho\log\E\exp 
\Big\{\rho^{-1}\theta\sqrt{\xi_1(\rho)+\cdots +\xi_l(\rho)}\Big\}\\ 
&&\le \sup_{\lambda >0}\Big\{\theta\sqrt{\lambda}-{\lambda\over 4 
C_2}\Big\} =C_2\theta^2.\nn 
\end{eqnarray} 
\qed

{\bf Proof of Theorem \ref{TLU1}:} 
\medskip We begin with the decomposition 
\begin{align} 
\big\vert Z[0,t]\big\vert &=\sum_{k=1}^{2^N}\bigg\vert 
Z\Big[{k-1\over 2^N}t,\hskip.05in{k\over 2^N}t\Big]\bigg\vert\nn\\ & 
\hspace{.6in}~~~-\sum_{j=1}^N\sum_{k=1}^{2^{j-1}}\bigg\vert 
Z\Big[{2k-2\over 2^j}t,\hskip.05in{2k-1\over 2^j}t\Big]\cap Z\Big[{2k-1\over 
2^j}t,\hskip.05in{2k\over 2^j}t\Big] 
\bigg\vert\label{4.102}\\ &=: I_t-J_t.\nn 
\end{align} 
 
We first establish the upper bound. Let $\eps>0$ be fixed. Since 
\begin{equation} E\big\vert Z[0,t]\big\vert -\big\vert Z[0,t]\big\vert =(\E 
I_t-I_t)+J_t -\E J_t\le (\E I_t-I_t)+J_t, 
\label{4.121} 
\end{equation} 
it follows that 
\begin{eqnarray} && 
\P\Big\{\Big\vert \E\big\vert Z[0,t]\big\vert -\big\vert Z[0,t] 
\big\vert\Big\vert\ge{\lambda t b_t/ (\log t)^2} 
\Big\}\label{4.140}\\ &&~~~~~~\le 
\P\Big\{\big\vert \E I_t -I_t\big\vert \ge{\eps 
  t b_t/ (\log t)^2} 
\Big\} +\P\Big\{J_t \ge{(\lambda-\eps) t b_t/ (\log t)^2} 
\Big\}.\nn 
\end{eqnarray} Notice that 
\begin{equation} 
\big\vert \E I_t -I_t\big\vert\le\sum_{k=1}^{2^N}\Bigg\vert 
\E\bigg\vert Z\Big[{k-1\over 2^N}t,\hskip.05in{k\over 
2^N}t\Big]\bigg\vert -\bigg\vert Z\Big[{k-1\over 
2^N}t,\hskip.05in{k\over 2^N}t\Big]\bigg\vert\Bigg\vert. 
\label{4.122} 
\end{equation} 
 
  Replacing $t$ by $2^{-N}t$, $\la$ by $2^{N}\la$ and $b_t$ by $\wt 
b_t=:b_{ 2^{N}t}$ in (\ref{4.123}) we obtain 
\begin{equation} 
\limsup_{t\to\infty}{1\over b_t}\log\P\Big\{\Big\vert 
\E\big\vert Z[0,2^{-N}t]\big\vert -\big\vert 
Z[0,2^{-N}t]\big\vert\Big\vert\ge 
\lambda {tb_t\over (\log t)^2}\Big\} 
\le -2^NC\lambda. 
\label{4.124} 
\end{equation} 
  Hence by Lemma \ref{LU1}, 
\begin{equation} 
\limsup_{t\to\infty}{1\over b_t}\log 
\P\Big\{\big\vert \E I_t -I_t\big\vert \ge{\eps t b_t\over (\log t)^2} 
\Big\}\le -\eps C2^N. 
\label{4.125} 
\end{equation} 
 
By the triangle inequality, 
\begin{equation} 
\P\Big\{J_t \ge{(\lambda-\eps) 
  t b_t\over (\log t)^2} 
\Big\} 
\le\sum_{j=1}^N\P\bigg\{\sum_{k=1}^{2^{j-1}}\xi_{j,k}\ge 2^{-j} 
{(\lambda-\eps) 
  t b_t\over (\log t)^2}\bigg\} 
\label{4.126} 
\end{equation} 
  where for each $1\le j\le N$, 
\begin{equation} 
\xi_{j,k}(t)=\bigg\vert Z\Big[{2k-2\over 2^j}t,\hskip.05in{2k-1\over 
2^j}t\Big]\cap Z\Big[{2k-1\over 2^j}t,\hskip.05in{2k\over 2^j}t\Big] 
\bigg\vert, \qquad k=1,\cdots, 2^{j-1}, 
\label{4.127} 
\end{equation} 
  forms an i.i.d. sequence with the same distribution as 
\begin{equation} 
\vert Z[0, 2^{-j}t]\cap Z'[0,2^{-j}t]\vert. 
\label{4.128} 
\end{equation} 
  By Theorem 1 in \cite{Cb} (with $2^{-j}t$ instead of $t$), for any 
$\lambda >0$, 
\begin{align} 
\lim_{t\to\infty}&{1\over b_t}\log \P\Big\{ 
\vert Z[0, 2^{-j}t]\cap Z'[0,2^{-j}t]\vert\ge{\lambda t b_t\over (\log 
t)^2} 
\Big\}\nn\\ 
&=-2^j(2\pi)^{-2}\det(\Gamma)^{-1/2}\kappa(2,2)^{-4}\lambda. 
\label{4.129} 
\end{align} 
  Therefore, by Lemma \ref{LU1}, 
\begin{equation} 
\lim_{t\to\infty}{1\over b_t}\log 
\P\bigg\{\sum_{k=1}^{2^{j-1}}\xi_{j,k}\ge  {\lambda 
  t b_t\over (\log t)^2}\bigg\}=-2^j 
(2\pi)^{-2}\det(\Gamma)^{-1/2}\kappa(2,2)^{-4}\lambda. 
\label{4.130} 
\end{equation} 
   In particular, 
\begin{align} 
\lim_{t\to\infty}&{1\over b_t}\log 
\P\bigg\{\sum_{k=1}^{2^{j-1}}\xi_{j,k}\ge 2^{-j} {(\lambda-\eps) 
  t b_t\over (\log t)^2}\bigg\}\nn\\ 
&=-(2\pi)^{-2}\det(\Gamma)^{-1/2}\kappa(2,2)^{-4}\,(\lambda -\eps) 
\label{4.131} 
\end{align} 
and therefore by (\ref{4.126}) 
\begin{equation} 
\lim_{t\to\infty}{1\over b_t}\log 
\P\Big\{J_t \ge{(\lambda-\eps) 
  t b_t\over (\log t)^2} 
\Big\} 
=-(2\pi)^{-2}\det(\Gamma)^{-1/2}\kappa(2,2)^{-4}\,(\lambda -\eps). 
\label{4.130c} 
\end{equation} 
Combining (\ref{4.140}), (\ref{4.125}) and (\ref{4.130c}) and letting 
$\eps\rar 0$ we obtain 
\begin{align} 
\limsup_{t\to\infty}&{1\over b_t}\log 
\P\Big\{\Big\vert \E\big\vert Z[0,t]\big\vert -\big\vert Z[0,t] 
\big\vert\Big\vert\ge{\lambda t b_t\over (\log t)^2} 
\Big\}\nn\\ 
&\le -(2\pi)^{-2}\det(\Gamma)^{-1/2}\kappa(2,2)^{-4}\lambda . 
\label{4.132} 
\end{align} 
 
  By Varadhan's integral lemma \cite[Section 4.3] {DZ} 
\begin{eqnarray}&\hspace{ .3in} 
\limsup_{t\to\infty}{1\over b_t}\log \E\exp\bigg\{ 
\theta\sqrt{b_t\over t}(\log t) 
\Big\vert \E\big\vert Z[0,t]\big\vert -\big\vert Z[0,t] 
\big\vert\Big\vert^{1/2}\bigg\}\label{4.141}\\ &\le\sup_{\lambda 
>0}\bigg\{ 
\theta\lambda^{1/2}-(2\pi)^{-2}\det(\Gamma)^{-1/2}\kappa(2,2)^{-4}\lambda 
\bigg\}\nn\\ &=(\theta\pi)^2 
\sqrt{\det(\Gamma)}\kappa(2,2)^4.\nn 
\end{eqnarray} 
  (The uniform exponential integrability is provided by 
Lemma \ref{LU3}.) 
 
\medskip We now prove the lower bound. 
Using induction on  $N$, one can see that 
\begin{eqnarray}&&A_t(\eps)=:\Lambda_\eps \Big({t\over 
b_t}\Big)^{-2} 
\sum_{x\in \Z^2}\bigg[\sum_{y\in Z[0,t]} h_\eps\Big(\sqrt{b_t\over 
t}(x-y)\Big)\bigg]^2\nn\\ &&\le\Lambda_\eps 
\Big({t\over b_t}\Big)^{-2} 
\sum_{k=1}^{2^N}\sum_{x\in \Z^2}\bigg[\sum_{y\in 
Z\big[{k-1\over 
2^N}t,\hskip.05in{k\over 2^N}t\big]} 
h_\eps\Big(\sqrt{b_t\over t}(x-y)\Big)\bigg]^2\nn\\ 
&&+2\Lambda_\eps \Big({t\over b_t}\Big)^{-2} 
\sum_{j=1}^N\sum_{k=1}^{2^{j-1}}\sum_{x\in 
\Z^2}\bigg[\sum_{y\in  Z\big[{2k-2\over 2^j}t,\hskip.05in{2k-1\over 
2^j}t\big]} h_\eps\Big(\sqrt{b_t\over 
t}(x-y)\Big)\bigg]\nn\\ &&\hskip 1.6in\times \bigg[ 
\sum_{y'\in Z\Big[{2k-1\over 
2^j}t,\hskip.05in{2k\over 2^j}t\Big]} 
h_\eps\Big(\sqrt{b_t\over t}(x-y')\Big)\bigg]\nn\\ 
&&=:I_t(\eps)+2J_t(\eps).\label{4.142} 
\end{eqnarray} 
  Therefore, with $   I_t , J_t $ given by (\ref{4.102}) 
\begin{eqnarray}& 
\E\big\vert Z[0,t]\big\vert -\big\vert Z[0,t]\big\vert =(\E I_t-I_t)+J_t 
-\E J_t\label{4.143}\\ &\ge (\E I_t-I_t)+J_t(\eps)-\vert J_t 
-J_t(\eps)\vert -\E J_t\nn\\ &\hspace{ .8in}\geq  (\E I_t-I_t) -{1\over 2} 
I_t(\eps)-\vert J_t -J_t(\eps)\vert -\E J_t  +{1\over 2}A_t(\eps). \nn 
\end{eqnarray} 
We will see that the dominant contribution to the lower bound comes from 
$A_t(\eps)$. 
By the last display we see that 
\begin{equation} 
{1\over 2}A_t(\eps)\leq \big\vert\E\big\vert Z[0,t]\big\vert -\big\vert 
Z[0,t]\big\vert\big\vert+\vert\E I_t-I_t\vert+ 
{1\over 2} 
I_t(\eps)+\vert J_t -J_t(\eps)\vert +\E J_t.\label{4.143a} 
\end{equation} 
and consequently 
\begin{eqnarray}&& 
\bigg\vert{1\over 2}A_t(\eps)\bigg\vert^{ 1/2}\leq \bigg\vert\E\big\vert 
Z[0,t]\big\vert -\big\vert Z[0,t]\big\vert\bigg\vert^{ 1/2}+\vert\E 
I_t-I_t\vert^{ 1/2}\nn\\ 
&&\hspace{ 1in}+  \bigg\vert{1\over 2} 
I_t(\eps)\bigg\vert^{ 1/2}+\vert J_t -J_t(\eps)\vert^{ 1/2} +\vert\E 
J_t\vert^{ 1/2}.\label{4.143b} 
\end{eqnarray} 
  Notice that it follows from (\ref{4.23}) that  
\begin{equation} 
\E J_t\leq C_{ N}{t\over (\log t)^2}. 
\label{4.160} 
\end{equation} 
If $\bar p$  is such that $p^{ -1}+\bar p ^{ -1}=1 $, then by the generalized 
H\"older inequality with $f=\theta\sqrt{b_t\over t}\log t$ we have 
\begin{eqnarray} &&\bigg\|  \exp {f \over p}\big\vert{1\over 
2}A_t(\eps)\big\vert^{ 1/2}     \bigg\|_{ 1} 
\nn\\ && 
\leq e^{C_{ N}\sqrt{b_t}}\,\bigg\|  \exp {f \over p}\big\vert\E\big\vert 
Z[0,t]\big\vert -\big\vert Z[0,t]\big\vert\big\vert^{ 1/2}     \bigg\|_{ p} 
\,\cdot\bigg\|  \exp {f \over p}\vert\E 
I_t-I_t\vert^{ 1/2}     \bigg\|_{3\bar p}\nn\\ 
&&\hspace{ 1in} 
\cdot \bigg\|  \exp {f \over p}\big\vert{1\over 2} 
I_t(\eps)\big\vert^{ 1/2}    \bigg\|_{3\bar p} 
\,\cdot\bigg\|  \exp {f \over p}\vert J_t -J_t(\eps)\vert^{ 1/2}   \bigg  \|_{ 
3\bar p}\label{4.143c} 
\end{eqnarray} 
Taking the $p$-th power and noting that $\bar p/p=1/( p-1)$, this can be 
rewritten as 
\begin{align} 
\E\exp&\bigg\{\theta\sqrt{b_t\over t}(\log t) 
\Big\vert 
\E\big\vert Z[0,t]\big\vert -\big\vert 
Z[0,t]\big\vert\Big\vert^{1/2}\bigg\}\label{4.144}\\ &\ge e^{-C_{ N}\sqrt{b_t}} 
\bigg[\E\exp\Big\{{3\theta\over p-1}\sqrt{b_t\over t}(\log t) 
\vert \E I_t-I_t\vert^{1/2}\Big\}\bigg]^{-{p-1\over 3}}\nn\\ &\quad 
\times\bigg[\E\exp\Big\{{3\theta\over p-1}\sqrt{b_t\over t}(\log t) 
I_t(\eps)^{1/2}\Big\}\bigg]^{-{p-1\over 3}}\nn\\ &\quad 
\times\bigg[\E\exp\Big\{{3\theta\over p-1}\sqrt{b_t\over t}(\log t) 
\vert J_t-J_t(\eps)\vert^{1/2}\Big\}\bigg]^{-{p-1\over 3}}\nn\\ &\quad 
\times\Bigg[\E\exp\bigg\{{\theta\over 2p}\sqrt{b_t\over t} 
(\log t) 
\vert A_{ t}( \eps)\vert^{1/2}\bigg\}\Bigg]^p.\nn 
\end{align} 
 
By Lemma \ref{LU6} 
\begin{align} \lim_{t\to\infty}&{1\over b_t}\log 
\E\exp\bigg\{{\theta\over 2p}\sqrt{b_t\over t} 
(\log t) 
\vert A_{ t}( \eps)\vert^{1/2}\bigg\}\label{4.148}\\ & =\sup_{g\in 
{\cal F}}\bigg\{{\pi\theta\over p}\sqrt{\det(\Gamma)} 
\bigg(\int_{\R^2}\vert (g^2\ast h_\eps)(x)\vert^2dx\bigg)^{1/2} 
\nn\\ 
&~~~-{1\over 2}\int_{\R^2}\langle\nabla g(x),\Gamma\nabla g(x)\rangle 
dx\bigg\}.\nn 
\end{align} 
This will give the main contribution to (\ref{4.144}). We now bound the other 
factors in (\ref{4.144}). 
 
Using Lemma \ref{LU1} together with (\ref{4.141}) (with $t$ replaced by 
$2^{-N}t$, 
$\theta$ by 
$2^{-N/2}\theta$, and 
$b_t$ by $\wt 
b_t=:b_{ 2^{N}t}$) we can 
prove that for any 
$\theta >0$, 
\begin{equation} 
\limsup_{t\to\infty}{1\over b_t}\log \E\exp\Big\{\theta 
\sqrt{b_t\over t}(\log t) 
\vert \E I_t-I_t\vert^{1/2}\Big\}\le 2^{-N}C\theta^2. 
\label{4.161} 
\end{equation} 
 
Using (\ref{4.145}) and  Lemma \ref{LU1}, we see that 
\begin{equation} 
\limsup_{t\to\infty}{1\over b_t}\log 
\E\exp\Big\{\theta\sqrt{b_t\over t}(\log t) I_t(\eps)^{1/2}\Big\}\le 
2^{-N}C\theta^2, 
\label{4.162} 
\end{equation} 
  where $C>0$ does not depend on $\eps$. 
\medskip Notice that 
\begin{equation} 
\vert J_t-J_t(\eps)\vert\le 
\sum_{j=1}^N\sum_{k=1}^{2^{j-1}}\vert K_{j,k}(\eps)\vert, 
\label{4.163} 
\end{equation} 
  where 
\begin{align}K_{j,k}(\eps)&= 
\bigg\vert Z\Big[{2k-2\over 2^j}t,\hskip.05in{2k-1\over 2^j}t\Big]\cap 
Z\Big[{2k-1\over 2^j}t,\hskip.05in{2k\over 2^j}t\Big] 
\bigg\vert\nn\\ &\hspace{ 1in} -\Lambda_\eps \Big({t\over 
b_t}\Big)^{-2} 
\sum_{x\in \Z^2}\bigg[\sum_{y\in Z\big[{2k-2\over 2^j}t,\hskip.05in{2k-1\over 
2^j}t\big]} h_\eps\Big(\sqrt{b_t\over 
t}(x-y)\Big)\bigg]\nn\\ &\hskip2.2in\times \bigg[\sum_{y'\in 
Z\big[{2k-1\over 2^j}t,\hskip.05in{2k\over 2^j}t\big]} 
h_\eps\Big(\sqrt{b_t\over t}(x-y')\Big)\bigg].\label{4.146} 
\end{align} 
  For each $1\le j\le N$, $K_{j,1}(\eps),\cdots, K_{j, 2^{N-1}}(\eps)$ 
forms an i.i.d sequence with the same  distribution as $B^{ ( j)}_{ t}(\eps)$. 
It then follows from Lemma \ref{LU7} and  H\"older's inequality that 
\begin{equation} 
\limsup_{\eps\to 0}\limsup_{t\to\infty}{1\over b_t}\log 
\E\exp\Big\{\theta\sqrt{b_t\over t}(\log t) 
\vert J_t-J_t(\eps)\vert^{1/2}\Big\}=0. 
\label{4.165} 
\end{equation} 
 
  Hence 
\begin{align} 
\liminf_{t\to\infty}&{1\over b_t}\log \E\exp\bigg\{ 
\theta\sqrt{b_t\over t}(\log t) 
\Big\vert \E\big\vert Z[0,t]\big\vert -\big\vert Z[0,t] 
\big\vert\Big\vert^{1/2}\bigg\}\nn\\ &\ge -2^{-N+1}C{p-1\over 
3}\Big({3\theta\over p-1}\Big)^2\nn\\ &\quad -{p-1\over 
3}\limsup_{t\to\infty}{1\over b_t}\log 
\E\exp\Big\{{3\theta\over p-1}\sqrt{b_t\over t}\log t 
\vert J_t-J_t(\eps)\vert^{1/2}\Big\}\nn\\ &\quad +p\sup_{g\in {\cal 
F}}\bigg\{{\pi\theta\over p}\sqrt{\det(\Gamma)} 
\bigg(\int_{\R^2}\vert (g^2\ast h_\eps)(x)\vert^2dx\bigg)^{1/2} 
\nn\\ 
&\quad -{1\over 2}\int_{\R^2}\langle\nabla g(x),\Gamma\nabla g(x)\rangle 
dx\bigg\}.\nn 
\end{align} 
  Take limits on the right hand side in the following order: let 
$\eps\to 0^+$, (using (\ref{4.165})), $N\to\infty$, and then $p\to 1^+$. We 
obtain 
\begin{align} 
\liminf_{t\to\infty}&{1\over b_t}\log \E\exp\bigg\{ 
\theta\sqrt{b_t\over t}(\log t) 
\Big\vert \E\big\vert Z[0,t]\big\vert -\big\vert Z[0,t] 
\big\vert\Big\vert^{1/2}\bigg\}\label{4.150}\\ &\ge \sup_{g\in {\cal 
F}}\bigg\{\pi\theta\sqrt{\det(\Gamma)} 
\bigg(\int_{\R^2}\vert g(x)\vert^4dx\bigg)^{1/2} -{1\over 
2}\int_{\R^2}\langle\nabla g(x),\Gamma\nabla g(x)\rangle dx\bigg\}\nn\\ 
&=(\pi\theta)^2\sqrt{\det(\Gamma)}\sup_{f\in {\cal F}} 
\bigg\{ 
\bigg(\int_{\R^2}\vert f(x)\vert^4dx\bigg)^{1/2} -{1\over 
2}\int_{\R^2}\vert\nabla f(x)\vert^2 dx\bigg\}\nn\\ 
&=(\pi\theta)^2\sqrt{\det(\Gamma)}\kappa (2,2)^4,\nn 
\end{align} 
  where the second step follows from the substitution 
$g(x)=\sqrt{\vert\det (A)\vert}f(Ax)$ with the $2\times 2$ matrix 
$A$ satisfying 
\begin{equation} A^\tau \Gamma 
A=(\pi\theta)^2\sqrt{\det(\Gamma)}I_{2\times 2} 
\label{4.166} 
\end{equation} 
($I_{2\times 2}$ is the $2\times 2$ identity matrix), and where the 
last step follows from Lemma A.2 in \cite{C}. 
\qed

\section{Exponential asymptotics for the smoothed 
range}\label{sec-smoothrange} 
 
In order to prove Lemma \ref{LU6} we first obtain a weak convergence 
result. 
 
Let $\beta >0$ and write 
\begin{equation} 
A_{t,\bb} (\eps)=:\Lambda_\eps (t)^{-2}\sum_{x\in \Z^2} 
\bigg[\sum_{y\in Z[0,\beta t]} 
h_\eps\Big({x-y\over \sqrt{t}}\Big)\bigg]^2\label{5.01} 
\end{equation} 
and 
\begin{equation} 
B_{t,\bb} (\eps)=:\Lambda_\eps (t)^{-2}\sum_{x\in 
\Z^2}\bigg[\sum_{y\in Z[0,\beta t]} 
h_\eps\Big({x-y\over \sqrt{t}}\Big)\bigg]\bigg[\sum_{y'\in Z'[0,\beta t]} 
h_\eps\Big({x-y'\over \sqrt{t}}\Big)\bigg].\label{5.01b} 
\end{equation}

\medskip Let $W(t), W'(t)$ be independent planar Brownian motions, each 
  with covariance matrix 
$\Gamma$ and write 
\begin{equation} 
\alpha_\eps ([0,t]^2)=\int_0^t\!\!\int_0^t(h_\eps\ast h_\eps) 
(W(s)-W'(r)\big)dr\,  ds 
\label{4.41} 
\end{equation}and 
\begin{equation} 
\alpha ([0,t]^2)=\lim_{\eps\to 0} \alpha_\eps([0,t]^2). 
\label{4.42} 
\end{equation}

\begin{lemma}\label{LU5} 
\begin{eqnarray}&&{(\log t)^2\over t}\Big[\big\vert Z[0, \beta t]\cap Z'[0, 
\beta t]\big\vert -B_{t,\bb} (\eps) 
\Big]\label{4.43}\\ &&~~~~\buildrel d\over\longrightarrow 
(2\pi)^2\det(\Gamma)\Big[\alpha ([0,\bb ]^2)-\alpha_\eps 
([0, \bb ]^2)\Big]\nn 
\end{eqnarray} and 
\begin{equation} {(\log t)^2\over t}A_{ t,\bb}( \eps) 
\buildrel d\over \longrightarrow (2\pi)^2\det(\Gamma) 
\int_{\R^2}\bigg(\int_0^\beta h_\eps \big(W(s)-x)ds\bigg)^2dx. 
\label{4.44} 
\end{equation} as $t\to\infty$. 
\end{lemma} 
 
\proof To prove (\ref{4.43}), we consider the following result given on 
p.697 of \cite{LR}: if $Z^{ ( t)}( s)=:{Z(ts) \over \sqrt{t}}$ then 
\begin{align} 
\Big(Z^{ ( t)}(\cdot),& (Z')^{ ( t)}(\cdot), {(\log t)^2\over  
t}\big\vert Z[0, \beta 
t]\cap Z'[0, \beta t]\big\vert\Big)\nn\\ 
&\quad \buildrel d\over \longrightarrow \Big(W(\cdot),W'(\cdot), 
(2\pi)^2\det(\Gamma)\alpha ([0,\beta]^2)\Big) 
\label{4.45} 
\end{align} 
in the Skorohod topology as $t\to\infty$. Actually, the proof  in 
\cite{LR} is for the discrete time random walk, but a similar proof works for 
$Z$. 
 
Let $M>0$ be fixed for a moment. Notice that 
\begin{equation} p_{t,\eps}(x)\equiv 
\Lambda_\eps (t)^{-1}h_\eps\Big({x\over\sqrt{t}}\Big), 
\hskip.2in x\in\Z^2, 
\label{4.46} 
\end{equation} defines a probability density on $\Z^2$ and that 
\begin{equation} 
\widehat{p}_{t,\eps}\Big({\lambda\over\sqrt{t}}\Big) 
=\Lambda_\eps(t)^{-1}\sum_{x\in\Z^2}h_\eps \Big({x\over\sqrt{t}}\Big) 
\exp\Big\{i\lambda\cdot {x\over\sqrt{t}}\Big\} 
\longrightarrow\widehat{h}_\eps(\lambda) 
\label{4.47} 
\end{equation} uniformly on $[-M,M]^2$ as $t\to\infty$. Consequently  
the family 
\begin{equation} 
\psi_t(x, y)=\int_{[-M,M]^2}\Big\vert 
\widehat{p}_{t,\eps}\Big({\lambda\over\sqrt{t}}\Big)\Big\vert^2 
\bigg[\int_0^{\beta }e^{i\lambda\cdot x(s)}ds\bigg] 
\bigg[\int_0^{\beta }e^{-i\lambda\cdot y(s')}ds'\bigg]\,d\lambda 
\label{4.48} 
\end{equation} are convergent continuous functionals on 
$D\Big([0,\beta], \R^2\Big)\otimes D\Big([0,\beta], \R^2\Big)$. 
Therefore 
\begin{eqnarray}&&\bigg({1\over t^2}\int_{[-M,M]^2}\Big\vert 
\widehat{p}_{t,\eps}\Big({\lambda\over\sqrt{t}}\Big)\Big\vert^2\nn\\ && 
\hspace{ .5in}\bigg[\int_0^{\beta t 
}\exp\Big\{i\lambda\cdot{Z(s)\over\sqrt{t}}\Big\}ds\bigg] 
\bigg[\int_0^{\beta t 
}\exp\Big\{-i\lambda\cdot{Z'(s')\over\sqrt{t}}\Big\}ds' 
\bigg]\,d\lambda, \nn\\ &&\hskip3in {(\log t)^2\over t}\big\vert Z[0, 
\beta t]\cap Z'[0, 
\beta t]\big\vert\bigg)\nn\\ 
&&=\bigg(\int_{[-M,M]^2}\Big\vert 
\widehat{p}_{t,\eps}\Big({\lambda\over\sqrt{t}}\Big)\Big\vert^2\nn\\ && 
\hspace{ .5in}\bigg[\int_0^{\beta 
}\exp\Big\{i\lambda\cdot Z^{ ( t)}(s)\Big\}ds\bigg] 
\bigg[\int_0^{\beta 
}\exp\Big\{-i\lambda\cdot (Z')^{ ( t)}(s')\Big\}ds' 
\bigg]\,d\lambda, \nn\\ &&\hskip3in {(\log t)^2\over t}\big\vert Z[0, 
\beta t]\cap Z'[0, 
\beta t]\big\vert\bigg)\nn\\ 
&&\buildrel d\over \longrightarrow 
\bigg(\int_{[-M,M]^2} 
\vert \widehat{h}_\eps(\lambda)\vert^2 
\bigg[\int_0^{\beta }e^{i\lambda\cdot W(s)}ds\bigg] 
\bigg[\int_0^{\beta }e^{-i\lambda\cdot W'(s')}ds'\bigg]\,d\lambda, 
\nn\\ 
&&\hskip3in (2\pi)^2\det(\Gamma)\alpha ([0,\beta]^2)\bigg).\label{4.49} 
\end{eqnarray} 
 
Recall that by Lemma 3 in \cite{Cb}, 
\begin{equation} 
\sup_t\E\exp\Big\{\theta{\log t\over t}\vert 
Z[0,t]\vert\Big\}<\infty 
\label{4.71} 
\end{equation} 
for all $\theta >0$. 
We will show that uniformly in $\lambda\in [-M,M]^2$ 
\begin{align} 
\lim_{t\to\infty}&{1\over t^2}\E\bigg\vert\int_0^{\beta t } 
\exp\Big\{i\lambda\cdot{Z(s)\over\sqrt{t}}\Big\}ds  
\nn\\&-{\log t\over 
2\pi\sqrt{\det(\Gamma)}} 
\sum_{x\in Z[0,\beta t]}\exp\Big\{i\lambda\cdot 
{x\over\sqrt{t}}\Big\}\bigg\vert^2 =0. 
\label{4.50a} 
\end{align} 
  Using the inequality 
\[|AA'-BB'|\leq  |A( B-B')|+|( A-B)B'|,\] 
  the Cauchy-Schwarz inequality and (\ref{4.71}), we see from (\ref{4.50a}) 
  that  uniformly in $\lambda\in [-M,M]^2$ 
\begin{align}  
\lim_{t\to\infty}&{1\over t^2}\E\bigg\vert\bigg[\int_0^{\beta t 
}\exp\Big\{i\lambda\cdot{Z(s)\over\sqrt{t}}\Big\}ds\bigg] 
\bigg[\int_0^{\beta t 
}\exp\Big\{-i\lambda\cdot{Z'(s')\over\sqrt{t}}\Big\}ds' 
\bigg]\nn\\ & 
-\({\log t\over 
2\pi\sqrt{\det(\Gamma)}}\)^{ 2} 
\bigg[\sum_{x\in Z[0,\beta t]}\exp\Big\{i\lambda\cdot 
{x\over\sqrt{t}}\Big\}\bigg] 
\label{4.50b}\\ 
&\hspace{2 in} \bigg[\sum_{x\in Z'[0,\beta t]}\exp\Big\{i\lambda\cdot 
{x'\over\sqrt{t}}\Big\}\bigg]\bigg\vert\nn\\ 
&=0. 
\nonumber 
\end{align} 
Together with (\ref{4.49}) this shows that 
\begin{eqnarray}&&\bigg(\Big({\log t\over 2\pi t}\Big)^2 
\int_{[-M,M]^2}\Big\vert 
\widehat{p}_{t,\eps}\Big({\lambda\over\sqrt{t}}\Big)\Big\vert^2 
  \nn\\ &&\hskip1in 
\bigg[\sum_{x\in Z[0,\beta t]}\exp\Big\{i\lambda\cdot 
{x\over\sqrt{t}}\Big\}\bigg]\bigg[ 
\sum_{x'\in Z'[0,\beta t]}\exp\Big\{-i\lambda\cdot 
{x'\over\sqrt{t}}\Big\}\bigg]\,d\lambda, \nn\\ &&\hskip3in {(\log 
t)^2\over t}\big\vert Z[0, \beta t]\cap Z'[0, \beta t]\big\vert\bigg)\nn\\ 
&&\buildrel d\over 
\longrightarrow \bigg(\det(\Gamma) 
\int_{[-M,M]^2}\vert \widehat{h}_\eps(\lambda)\vert^2 
\bigg[\int_0^{\beta }e^{i\lambda\cdot W(s)}ds\bigg] 
\bigg[\int_0^{\beta }e^{-i\lambda\cdot W'(s')}ds'\bigg]\,d\lambda,\nn\\ 
&&\hskip3in (2\pi)^2\det(\Gamma)\alpha ([0,\beta]^2)\bigg).\label{4.53} 
\end{eqnarray}

Notice by (\ref{4.47}) that for any $\delta >0$, one can take $M>0$ 
sufficiently large so that 
\begin{equation} 
\Big\vert 
\widehat{p}_{t,\eps}\Big({\lambda\over\sqrt{t}}\Big)\Big\vert <\delta, 
\hskip.2in\lambda\in [-\sqrt{t}\pi,\sqrt{t}\pi]^2\setminus[M,M]^2, 
\label{4.54} 
\end{equation} 
  if $t$ is sufficiently large. Consequently 
\begin{eqnarray}&&H_{ t}=: 
\bigg\vert\int_{[-\sqrt{t}\pi,\sqrt{t}\pi]^2\setminus[-M,M]^2} 
\Big\vert 
\widehat{p}_{t,\eps}\Big({\lambda\over\sqrt{t}}\Big)\Big\vert^2 
\bigg[\sum_{x\in Z[0,\beta t]}\exp\Big\{i\lambda\cdot 
{x\over\sqrt{t}}\Big\}\bigg]\label{4.54A}\\ 
&&\hskip2.1in\times\bigg[ 
\sum_{x'\in Z'[0,\beta t]}\exp\Big\{-i\lambda\cdot 
{x'\over\sqrt{t}}\Big\}\bigg]\,d\lambda\bigg\vert\nn\\ &&\le (2\pi)^2 
\delta t\big\vert Z[0, \beta t]\cap Z'[0, \beta t]\big\vert.\nn 
\end{eqnarray} 
It follows from  (\ref{4.23}) that $\Big({\log t/( 2\pi t)}\Big)^2H_{ t}\rar 
0$ in $L^{ 1}$ uniformly in large $t$ as $M\to\infty$.  
Therefore, using (\ref{4.53}) and the fact that 
$\widehat{h}\in L^{ 2}$, we obtain 
\begin{eqnarray}&&\bigg(\Big({\log t\over 2\pi t}\Big)^2 
\int_{[-\sqrt{t}\pi,\sqrt{t}\pi]^2}\Big\vert 
\widehat{p}_{t,\eps}\Big({\lambda\over\sqrt{t}}\Big)\Big\vert^2 
\bigg[\sum_{x\in Z[0,\beta t]}\exp\Big\{i\lambda\cdot 
{x\over\sqrt{t}}\Big\}\bigg]\nn\\ &&\hskip.3in\times\bigg[ 
\sum_{x'\in Z'[0,\beta t]}\exp\Big\{-i\lambda\cdot 
{x'\over\sqrt{t}}\Big\}\bigg]\,d\lambda, \hskip.1in {(\log t)^2\over 
t}\big\vert Z[0, 
\beta t]\cap Z'[0, \beta t]\big\vert\bigg)\nn\\ &&\buildrel d\over 
\longrightarrow \bigg(\det(\Gamma) 
\int_{\R^2} \vert \widehat{h}_\eps(\lambda)\vert^{ 2} 
\bigg[\int_0^{\beta }e^{i\lambda\cdot W(s)}ds\bigg] 
\bigg[\int_0^{\beta }e^{-i\lambda\cdot W'(s')}ds'\bigg]\,d\lambda,\nn\\ 
&&\hskip.6in (2\pi)^2\det(\Gamma)\alpha ([0,\beta]^2)\bigg).\label{4.55} 
\end{eqnarray} 
Note that 
\begin{eqnarray} && 
B_{t,\bb} (\eps)=\sum_{x\in 
\Z^2}\bigg[\sum_{y\in Z[0,\beta t]}\Lambda_\eps (t)^{-1} 
h_\eps\Big({x-y\over \sqrt{t}}\Big)\bigg]\bigg[\sum_{y'\in Z'[0,\beta t]} 
\Lambda_\eps (t)^{-1}h_\eps\Big({x-y'\over 
\sqrt{t}}\Big)\bigg]\nn\\ &&\hspace{ .5in} = \sum_{y\in Z[0,\beta 
t]}\sum_{y'\in Z'[0,\beta t]} \bigg[\sum_{x\in 
\Z^2} 
p_{t,\eps}(x-y) 
p_{t,\eps}(x-y')\bigg].\label{4.55h} 
\end{eqnarray} 
It then follows from Parseval's identity that  
\begin{eqnarray}&&(2\pi)^2t\,B_{t,\bb} (\eps).\label{4.56}\\ 
&&=t\int_{[-\pi,\pi]^2}  \vert 
\widehat{p}_{t,\eps}(\lambda)\vert^2 
\bigg[\sum_{y\in Z[0,\beta t]}e^{i\lambda\cdot y}\bigg] 
\bigg[\sum_{y'\in Z'[0,\beta t]}e^{-i\lambda\cdot y'}\bigg]\,d\lambda\nn\\ 
&& =\int_{[-\sqrt{t}\pi,\sqrt{t}\pi]^2} \Big\vert 
\widehat{p}_{t,\eps}\Big({\lambda\over\sqrt{t}}\Big)\Big\vert^2 
\nn\\ 
&&\hspace{.5in} 
\bigg[\sum_{y\in Z[0,\beta t]}\exp\Big\{i\lambda\cdot 
{y\over\sqrt{t}}\Big\}\bigg]\bigg[\sum_{y'\in Z'[0,\beta t]} 
\exp\Big\{-i\lambda\cdot  {y'\over\sqrt{t}}\Big\}\bigg]\,d\lambda.\nn 
\end{eqnarray} 
Similarly, using the fact that $h_\eps$ is symmetric so that 
$\widehat{h}_\eps(\lambda)$ is real 
\begin{equation} 
\int_{\R^2} \vert \widehat{h}_\eps(\lambda)\vert^{ 2} 
\bigg[\int_0^{\beta }e^{i\lambda\cdot W(s)}ds\bigg] 
\bigg[\int_0^{\beta }e^{-i\lambda\cdot W'(s')}ds'\bigg]\,d\lambda 
=\alpha_\eps ([0,\beta]^2).\label{4.56A} 
\end{equation} 
Thus, we have proved 
\begin{eqnarray} &&\bigg({(\log t)^2\over t}B_{t,\bb} (\eps),\hskip.1in 
{(\log t)^2\over t}\big\vert Z[0, \beta t]\cap Z'[0, \beta 
t]\big\vert\bigg)\label{4.58}\\ &&\buildrel d\over \longrightarrow 
\Big((2\pi)^2\det(\Gamma)\alpha_\eps ([0,\beta]^2),\hskip.1in 
(2\pi)^2\det(\Gamma)\alpha([0,\beta]^2)\Big).\nn 
\end{eqnarray} 
  (\ref{4.43}) follows from this.

Thus to complete the proof of  (\ref{4.43}) it only remains to show 
(\ref{4.50a}) uniformly in 
$\lambda\in [-M,M]^2$. 
We will show that for  any $\de>0$ we can find $\de'>0$ and  $t_{ 0}<\ff$ 
such that 
\begin{equation} 
{1\over t^2}\E\bigg\vert\int_0^{\beta t} 
\exp\Big\{i\lambda\cdot {Z(s)\over\sqrt{t}}\Big\}ds-\int_0^{\beta t} 
\exp\Big\{i\gamma\cdot {Z(s)\over\sqrt{t}}\Big\}ds\bigg\vert^2 
<\delta\label{4.50x} 
\end{equation} 
and 
\begin{equation} 
\Big({\log t\over t}\Big)^2\E\bigg\vert\sum_{x\in Z[0,\beta t]} 
\exp\Big\{i\lambda\cdot {x\over\sqrt{t}}\Big\}-\sum_{x\in Z[0,\beta t]} 
\exp\Big\{i\gamma\cdot {x\over\sqrt{t}}\Big\}\bigg\vert^2 
<\delta\label{4.50y} 
\end{equation} 
for all $t\geq t_{ 0}$ and 
$\vert\lambda -\gamma\vert\leq \de'$. 
We then cover $[-M, M]^2$ by a finite number of discs $B(\lambda_k,\de' 
)$ of radius $\de' $ centered at $\lambda_k$, $k=1,\ldots,N$. Define 
$\tau (\lambda)=\lambda_k$ where $k$ is the smallest integer with 
$\lambda\in B(\lambda_k,\de' 
)$. By \cite[(4.11)]{Cb}, we can choose $t_{ 1}<\ff$ such that for all 
$t\geq t_{ 1}$ and $k=1,\ldots,N$, 
\begin{equation} 
{1\over t^2}\E\bigg\vert\int_0^{\beta t } 
\exp\Big\{i\lambda_k\cdot{Z(s)\over\sqrt{t}}\Big\}ds -{\log t\over 
2\pi\sqrt{\det(\Gamma)}} 
\sum_{x\in Z[0,\beta t]}\exp\Big\{i\lambda_k\cdot 
{x\over\sqrt{t}}\Big\}\bigg\vert^2 \leq \de. 
\label{4.50r} 
\end{equation} 
Hence, uniformly in  $\lambda\in [-M,M]^2$ we have that for all $t\geq t_{ 
0}\vee t_{ 1}$ 
\begin{equation} 
{1\over t^2}\E\bigg\vert\int_0^{\beta t } 
\exp\Big\{i\lambda\cdot{Z(s)\over\sqrt{t}}\Big\}ds -{\log t\over 
2\pi\sqrt{\det(\Gamma)}} 
\sum_{x\in Z[0,\beta t]}\exp\Big\{i\lambda\cdot 
{x\over\sqrt{t}}\Big\}\bigg\vert^2 \leq 3\de 
\label{4.50s} 
\end{equation} 
proving that (\ref{4.50a}) holds uniformly in 
$\lambda\in [-M,M]^2$.

  (\ref{4.50x}) actually holds uniformly in $t$. To see this note that 
\begin{eqnarray} && 
{1\over t^2}\E\bigg\vert\int_0^{\beta t } 
\exp\Big\{i\lambda\cdot{Z(s)\over\sqrt{t}}\Big\}ds- 
\int_0^{\beta t } 
\exp\Big\{i\gamma\cdot{Z(s)\over\sqrt{t}}\Big\}ds\bigg\vert^2 
\label{4.51}\\ && 
\leq  {1\over t^2}\E\bigg\vert\int_0^{\beta t } 
|\lambda-\gamma|{|Z(s)|\over\sqrt{t}}ds\bigg\vert^2\nonumber\\ && 
=  {|\lambda-\gamma|^{ 2}\over t^3}\E 
\int_0^{\beta t }\int_0^{\beta t } |Z(s)| |Z(r)|\,ds\,dr\nonumber\\ && 
\leq  C{|\lambda-\gamma|^{ 2}\over t^3} 
\int_0^{\beta t }\int_0^{\beta t }s^{ 1/2}r^{ 1/2}\,ds\,dr\leq C' 
|\lambda-\gamma|^{ 2}.\nonumber 
\end{eqnarray} 
As for (\ref{4.50y}), 
\begin{eqnarray} && 
\E\bigg\vert\sum_{x\in Z[0,\beta t]} 
\exp\Big\{i\lambda\cdot {x\over\sqrt{t}}\Big\}-\sum_{x\in Z[0,\beta t]} 
\exp\Big\{i\gamma\cdot {x\over\sqrt{t}}\Big\}\bigg\vert^2\label{4.51y}\\ 
&&\le  4\E \lc\vert Z[0,\beta t]\vert^2 
1_{\{\sup_{s\le\beta t}\vert Z(s)\vert\ge C\sqrt{t}\}}\rc\nonumber\\ 
&&+|\lambda-\gamma|^{ 2} 
\E\lc\Big\vert  \sum_{x\in Z[0,\beta t]}{|x|\over\sqrt{t}} 
\Big\vert^2 
1_{\{\sup_{s\le\beta t}\vert Z(s)\vert\le C\sqrt{t}\}}\rc\nonumber\\ 
&&\le 4\E \lc\vert Z[0,\beta t]\vert^2 
1_{\{\sup_{s\le\beta t}\vert Z(s)\vert\ge C\sqrt{t}\}}\rc 
+C^2\vert \lambda -\gamma\vert^2\E \vert Z[0,\beta t]\vert^2\nonumber 
\end{eqnarray} 
and by (\ref{4.71}) 
\begin{eqnarray} &&\quad 
4\E \lc\vert Z[0,\beta t]\vert^2 
1_{\{\sup_{s\le\beta t}\vert Z(s)\vert\ge C\sqrt{t}\}}\rc 
+C^2\vert \lambda -\gamma\vert^2\E \vert Z[0,\beta 
t]\vert^2\label{4.51z}\\ && \leq 
4\lc\E (\vert Z[0,\beta t]\vert^4) 
P(\sup_{s\le\beta t}\vert Z(s)\vert\ge C\sqrt{t})\rc^{ 1/2} 
+C^2\vert \lambda -\gamma\vert^2\E \vert Z[0,\beta 
t]\vert^2\nonumber\\ && \leq \({ ct\over \log t}\)^{ 2} 
\Bigg(4\lc 
P(\sup_{s\le\beta t}\vert Z(s)\vert\ge C\sqrt{t})\rc^{ 1/2} 
+C^2\vert \lambda -\gamma\vert^2\Bigg).\nonumber 
\end{eqnarray} 
Taking $C$ large and then choosing $\de'>0$ sufficiently small 
  completes the proof of (\ref{4.50y}) and hence of  (\ref{4.43}). 
\medskip

We now prove (\ref{4.44}). Using the facts that $\Lambda_\eps (t)\sim t$, 
that 
\begin{eqnarray} &&{1\over t}\sum_{x\in \Z^2}\bigg[\sum_{y\in Z[0,\beta 
t]} h_\eps\Big({x-y\over \sqrt{t}}\Big)\bigg]^2 
-\int_{\R^2}\bigg[\sum_{y\in Z[0,\beta t]} h_\eps\Big(x-{y\over 
\sqrt{t}}\Big)\bigg]^2dx\nn\\ 
&&\hspace{ 3in}=o(1)\vert Z[0,\beta t]\vert^2,\label{4.60} 
\end{eqnarray} (where the boundedness and continuity of $h_\eps$ is used), 
and (\ref{4.71}) we need only show that 
\begin{align} 
\Big({\log t\over t}\Big)^2&\int_{\R^2}\bigg[\sum_{y\in Z[0,\beta t]} 
h_\eps\Big(x-{y\over \sqrt{t}}\Big)\bigg]^2dx\\ 
&\buildrel d\over \longrightarrow (2\pi)^2\det(\Gamma) 
\int_{\R^2}\bigg(\int_0^\beta h_\eps \big(W(s)-x)ds\bigg)^2dx.\nn 
\end{align} 
  By the Parseval identity, 
\begin{align}  \int_{\R^2}&\bigg[\sum_{y\in Z[0,\beta t]} 
h_\eps\Big(x-{y\over 
\sqrt{t}}\Big)\bigg]^2dx\label{4.61}\\ 
&=(2\pi)^{-2}\int_{\R^2} \bigg\vert\int_{\R^2}e^{i\lambda\cdot 
x} 
\sum_{y\in Z[0,\beta t]} h_\eps\Big(x-{y\over \sqrt{t}}\Big)\,dx 
\bigg\vert^2\,d\lambda\nn\\ 
&=(2\pi)^{-2}\int_{\R^2} \bigg\vert\int_{\R^2} h_\eps (x) 
e^{i\lambda\cdot x}dx\bigg\vert^2\bigg\vert\sum_{y\in Z[0,\beta t]} 
\exp\Big\{i\lambda\cdot 
{y\over\sqrt{t}}\Big\}\bigg\vert^2\,d\lambda\nn\\ 
&=\int_{\R^2} \big\vert \wh h_\eps (\la) 
\big\vert^2\big\vert\sum_{y\in Z[0,\beta t]} 
\exp\Big\{i\lambda\cdot {y\over\sqrt{t}}\Big\}\big\vert^2\,d\lambda.\nn 
\end{align} 
 
Let $M>0$ be fixed and 
$\lambda_1,\cdots,\lambda_N$ and $\tau$ be defined as above. 
  By  \cite[Theorem 7]{Cb}, 
\begin{eqnarray} && 
{\log t\over t}\bigg(\sum_{y\in Z[0,\beta t]} 
\exp\Big\{i\lambda_1\cdot {y\over\sqrt{t}}\Big\},\cdots, 
\sum_{y\in Z[0,\beta t]} 
\exp\Big\{i\lambda_N\cdot {y\over\sqrt{t}}\Big\}\bigg)\label{4.61a}\\ 
&&~~~\buildrel d\over\longrightarrow (2\pi)\sqrt{\det(\Gamma)} 
\bigg(\int_0^\beta e^{i\lambda_1\cdot W(s)}ds,\cdots, 
\int_0^\beta e^{i\lambda_N\cdot W(s)}ds\bigg).\nn 
\end{eqnarray} 
In particular, 
\begin{eqnarray} && 
\Big({\log t\over t}\Big)^2\int_{[-M,M]^2}\vert \widehat{h}_\eps 
(\lambda)\vert^2 
\bigg\vert\sum_{y\in Z[0,\beta t]} 
\exp\Big\{i\tau(\lambda)\cdot {y\over\sqrt{t}}\Big\}\bigg\vert^2d 
\lambda\label{4.61b}\\ 
&&=\sum_{k=1}^N\int_{B_k}\vert \widehat{h}_\eps 
(\lambda)\vert^2 
\bigg\vert{\log t\over t}\sum_{y\in Z[0,\beta t]} 
\exp\Big\{i\lambda_k\cdot {y\over\sqrt{t}}\Big\}\bigg\vert^2d\lambda\nn\\ 
&&\buildrel d\over\longrightarrow (2\pi)^2\det(\Gamma) 
\sum_{k=1}^N\int_{B_k}\vert \widehat{h}_\eps (\lambda)\vert^2 
\bigg\vert\int_0^\beta e^{i\lambda_k\cdot W(s)}ds\bigg\vert^2d\lambda\nn\\ 
&&=(2\pi)^2\det(\Gamma)\int_{[-M,M]^2}\vert 
\widehat{h}_\eps (\lambda)\vert^2 
\bigg\vert\int_0^\beta e^{i\tau(\lambda)\cdot W(s)}ds 
\bigg\vert^2d\lambda.\nn 
\end{eqnarray} 
Notice that the right hand side of (\ref{4.61b}) converges to 
\begin{equation} 
(2\pi)^2\det(\Gamma)\int_{[-M,M]^2}\vert 
\widehat{h}_\eps(\lambda)\vert^2 
\bigg\vert\int_0^\beta e^{i\lambda\cdot 
W(s)}ds\bigg\vert^2d\lambda\label{4.61cc} 
\end{equation} 
as $N\to\infty$. Applying (\ref{4.50y}) to the left hand side of (\ref{4.61b}) 
gives 
\begin{eqnarray} && 
\Big({\log t\over t}\Big)^2\int_{[-M,M]^2}\vert \widehat{h}_\eps 
(\lambda)\vert^2 
\bigg\vert\sum_{y\in Z[0,\beta t]} 
\exp\Big\{i\lambda\cdot 
{y\over\sqrt{t}}\Big\}\bigg\vert^2d\lambda\label{4.61c}\\ &&\buildrel 
d\over\longrightarrow(2\pi)^2\det(\Gamma)\int_{[-M,M]^2}\vert 
\widehat{h}_\eps (\lambda)\vert^2 
\bigg\vert\int_0^\beta e^{i\lambda\cdot W(s)}ds\bigg\vert^2d\lambda.\nn 
\end{eqnarray} 
As $M\to\infty$, the right hand side of (\ref{4.61c}) converges to 
\begin{eqnarray} && 
(2\pi)^2\det(\Gamma)\int_{\R^2}\vert 
\widehat{h}_\eps (\lambda)\vert^2 
\bigg\vert\int_0^\beta e^{i\lambda\cdot 
W(s)}ds\bigg\vert^2d\lambda\label{4.61d}\\ 
&&=\det(\Gamma)\int_{\R^2}\bigg(\int_0^\beta h_\epsilon \big( 
W(s)-x\big)ds\bigg)^2dx \nn 
\end{eqnarray}
by Parseval's identity. 
Note 
\begin{eqnarray} && 
H'_{t,M}=:\int_{\R^2\setminus[-M,M]^2}\vert 
\widehat{h}_\eps (\lambda)\vert^2 
\bigg\vert\sum_{y\in Z[0,\beta t]} 
\exp\Big\{i\lambda\cdot 
{y\over\sqrt{t}}\Big\}\bigg\vert^2d\lambda\label{4.61e}\\ &&\le \big\vert 
Z[0,\beta t]\big\vert^2 
\int_{\R^2\setminus[-M,M]^2}\vert \widehat{h}_\eps 
(\lambda)\vert^2d\lambda.\nn 
\end{eqnarray} 
It follows from  (\ref{4.71}) and the fact that $\widehat{h}_\eps \in L^{ 2}$ 
that 
$\Big({\log t\over 2\pi t}\Big)^2H'_{t,M}\rar 0$ in $L^{ 1}$ as $M\rar \ff$ 
uniformly in $t$. Therefore, using the last three displays, we obtain 
\begin{eqnarray} && 
\Big({\log t\over t}\Big)^2\int_{\R^2}\vert \widehat{h}_\eps 
(\lambda)\vert^2 
\bigg\vert\sum_{y\in Z[0,\beta t]} 
\exp\Big\{i\lambda\cdot 
{y\over\sqrt{t}}\Big\}\bigg\vert^2d\lambda\label{4.61f}\\ &&\buildrel 
d\over\longrightarrow 
\det(\Gamma)\int_{\R^2}\bigg(\int_0^\beta h_\epsilon \big( 
W(s)-x\big)ds\bigg)^2dx.\nn 
\end{eqnarray} 
\qed

{\bf  Proof of Lemma \ref{LU6}}: Let $T>0$ be fixed for the moment. Write 
$\gamma_t=t/[T^{-1}b_t]$. We have 
\begin{eqnarray} &&\E\exp\bigg\{\theta\sqrt{b_t\over t} 
(\log t) 
\bigg(\Lambda_\eps\Big({t\over b_t}\Big)^{-2}\sum_{x\in 
\Z^2}\bigg[\sum_{y\in Z[0,t]} h_\eps\Big(\sqrt{b_t\over 
t}(x-y)\Big)\bigg]^2 
\bigg)^{1/2}\bigg\}\nn\\ &&\le \Bigg[\E\exp\bigg\{\theta\sqrt{b_t\over t} 
(\log t)\nn\\ &&\hskip1in\times\bigg( 
\Lambda_\eps\Big({t\over b_t}\Big)^{-2}\sum_{x\in 
\Z^2}\bigg[\sum_{y\in Z[0,\gamma_t]} h_\eps\Big(\sqrt{b_t\over 
t}(x-y)\Big)\bigg]^2 
\bigg)^{1/2}\bigg\}\Bigg]^{[T^{-1}b_t]}.\label{4.68} 
\end{eqnarray} 
We obtain from Lemma \ref{LU5} (with $t$ being replaced by $t/b_t$ and 
$\bb=T$) 
\begin{eqnarray} &&{b_t\over t}(\log t)^2\Lambda_\eps\Big({t\over 
b_t}\Big)^{-2} 
\sum_{x\in \Z^2}\bigg[\sum_{y\in Z[0,\gamma_t]} 
h_\eps\Big(\sqrt{b_t\over t}(x-y)\Big)\bigg]^2\nn 
\\ 
&&\hspace{ 1in}\buildrel d\over \longrightarrow (2\pi)^2\det(\Gamma) 
\int_{\R^2}\bigg(\int_0^T h_\eps \big(W(s)-x)ds\bigg)^2dx,\label{4.69} 
\hskip.2in t\to\infty. 
\end{eqnarray} 
In addition, 
\begin{align}{b_t\over t}&(\log t)^2\Lambda_\eps\Big({t\over 
b_t}\Big)^{-2}  \sum_{x\in \Z^2}\bigg[\sum_{y\in Z[0,\gamma_t]} 
h_\eps\Big(\sqrt{b_t\over t}(x-y)\Big)\bigg]^2\nn 
\\ &\le{b_t\over t}(\log t)^2\Lambda_\eps\Big({t\over b_t}\Big)^{-2} 
\vert\vert h\vert\vert_\infty 
\vert Z[0,\gamma_t]\vert\sum_{\stackrel{x\in\Z^2}{y\in Z[0,\gamma_t]}} 
h_\eps\Big(\sqrt{b_t\over t}(x-y)\Big)\nn\\ 
&={b_t\over t}(\log t)^2\Lambda_\eps\Big({t\over b_t}\Big)^{-2} 
\vert\vert h\vert\vert_\infty 
\vert Z[0,\gamma_t]\vert^2\sum_{x\in \Z^2} 
h_\eps\Big(\sqrt{b_t\over t}x\Big)\nn\\ &\le 
C\Big({b_t\over t}\Big)^2 (\log t)^2\vert Z[0,\gamma_t]\vert^2, \label{4.70} 
\end{align} 
where in the last step we used (\ref{4.133n}). 
  (\ref{4.71}) together with (\ref{4.69}) then implies that 
\begin{eqnarray} &\E\exp\bigg\{\theta\sqrt{b_t\over t} 
\Lambda_\eps\Big({t\over b_t}\Big)^{-1} 
(\log t) 
\bigg(\sum_{x\in \Z^2}\bigg[\sum_{y\in Z[0,\gamma_t]} 
h_\eps\Big(\sqrt{b_t\over t}(x-y)\Big)\bigg]^2 
\bigg)^{1/2}\bigg\}\nn\\ 
&\hspace{ 
.5in}\longrightarrow\E\exp\bigg\{2\pi\theta\sqrt{\det(\Gamma)} 
\bigg(\int_{\R^2}\bigg(\int_0^T h_\eps 
  \big(W(s)-x)ds\bigg)^2dx\bigg)^{1/2}\bigg\}.\label{4.72} 
\end{eqnarray} 
  Combining (\ref{4.68})  and (\ref{4.72}) we see that 
\begin{align} \limsup_{t\to\infty}&{1\over b_t}\log 
\E\exp\bigg\{\theta\sqrt{b_t\over t} 
(\log t) 
\vert A_{ t}(\eps )\vert^{1/2}\bigg\}\label{4.73}\\ &\le {1\over T}\log 
\E\exp\bigg\{2\pi\theta\sqrt{\det(\Gamma)} 
\bigg(\int_{\R^2}\bigg(\int_0^T h_\eps 
\big(W(s)-x)ds\bigg)^2dx\bigg)^{1/2} 
\bigg\}.\nn 
\end{align} 
  Then the upper bound for (\ref{1.00a}) follows from the fact  that 
\begin{align} \lim_{T\to\infty}&{1\over T}\log 
\E\exp\bigg\{2\pi\theta\sqrt{\det(\Gamma)} 
\bigg(\int_{\R^2}\bigg(\int_0^T h_\eps 
\big(W(s)-x)ds\bigg)^2dx\bigg)^{1/2} 
\bigg\}\nn\\ &=\sup_{g\in {\cal F}}\bigg\{2\pi\theta\sqrt{\det(\Gamma)} 
\bigg(\int_{\R^2}\vert (g^2\ast h_\eps)(x)\vert^2dx\bigg)^{1/2} 
\nn\\&\hspace{2 in} -{1\over 2}\int_{\R^2}\langle\nabla g(x), \Gamma\nabla g(x)\rangle 
dx\bigg\}.\nn 
\end{align} 
This is  \cite[Theorem 7]{CR}. (Or see the earlier 
\cite[Theorem 3.1]{C}, which uses a slightly different smoothing). 
 
We now prove the lower bound for (\ref{1.00a}). Let $f$ be a smooth 
function on $\R^2$ with compact support and 
\begin{equation} 
\vert\vert f\vert\vert_2=\bigg(\int_{\R^2}\vert 
f(x)\vert^2dx\bigg)^{1/2}=1. 
\label{4.74} 
\end{equation} 
We can write 
\begin{eqnarray} && 
\sqrt{b_t\over t}\bigg(\sum_{x\in \Z^2}\bigg[\sum_{y\in Z[0,t]} 
h_\eps\Big(\sqrt{b_t\over t}(x-y)\Big)\bigg]^2 
\bigg)^{1/2}\label{4.74a}\\ && 
=\sqrt{b_t\over t}\bigg(\int_{\R^2}\bigg[\sum_{y\in Z[0,t]} 
h_\eps\Big(\sqrt{b_t\over t}([x]-y)\Big)\bigg]^2 
\,dx\bigg)^{1/2}\nonumber\\ && 
=\bigg(\int_{\R^2}\bigg[\sum_{y\in Z[0,t]} 
h_\eps\Big(\sqrt{b_t\over t}(\Big[\sqrt{t\over b_t}x\Big]-y)\Big)\bigg]^2 
\,dx\bigg)^{1/2}.\nonumber 
\end{eqnarray} 
Hence by the  Cauchy-Schwarz inequality, 
\begin{align} 
\sqrt{b_t\over t}&\bigg(\sum_{x\in \Z^2}\bigg[\sum_{y\in Z[0,t]} 
h_\eps\Big(\sqrt{b_t\over t}(x-y)\Big)\bigg]^2 
\bigg)^{1/2}\nn\\ &=\bigg(\int_{\R^2}\bigg[\sum_{y\in Z[0,t]} 
h_\eps\Big(\sqrt{b_t\over t}\Big[\sqrt{t\over b_t}x\Big]-\sqrt{b_t\over 
t}y 
\Big)\bigg]^2dx 
\bigg)^{1/2}\nn\\ &\ge \int_{\R^2}f(x)\sum_{y\in Z[0,t]} 
h_\eps\Big(\sqrt{b_t\over t}\Big[\sqrt{t\over b_t}x\Big]-\sqrt{b_t\over 
t}y 
\Big)dx\nn\\ &=\int_{\R^2}f(x)\sum_{y\in Z[0,t]} 
h_\eps\Big(x-\sqrt{b_t\over t}y 
\Big)dx +O(1)\vert Z[0,t]\vert,\hskip.2in t\to\infty,\label{4.75} 
\end{align} 
where $O(1)$ is bounded by a constant. In view of (4.12),  
recalling that $$\sqrt{b_t\over t} 
\vert A_{ t}(\eps )\vert^{1/2}\sim {b_t\over t}\sqrt{b_t\over 
t}\bigg(\sum_{x\in \Z^2}\bigg[\sum_{y\in Z[0,t]} 
h_\eps\Big(\sqrt{b_t\over t}(x-y)\Big)\bigg]^2 
\bigg)^{1/2},$$ 
and using H\"older's 
inequality one can see that the term $O(1)\vert Z[0,t]\vert$ 
does not contribute anything to (\ref{1.00a}). 
 
\medskip

By  \cite[Theorem 8]{Cb}, 
\begin{align} \liminf_{t\to\infty}&{1\over b_t}\log 
\E\exp\bigg\{\theta{b_t\log t\over t} 
\sum_{y\in Z[0, t]}(f\ast h_\eps)\Big(\sqrt{b_t\over 
t}y\Big)\bigg\}\label{4.76}\\ &\ge\sup_{g\in {\cal 
F}}\bigg\{2\pi\theta\sqrt{\det(\Gamma)} 
\int_{\R^2}(f\ast h_\eps)(x) g^2(x)dx -{1\over 
2}\int_{\R^2}\langle\nabla g(x),\Gamma\nabla g(x)\rangle dx\bigg\}\nn 
\\ &=\sup_{g\in {\cal F}}\bigg\{2\pi\theta\sqrt{\det(\Gamma)} 
\int_{\R^2}f(x)(g^2\ast h_\eps)(x) dx -{1\over 
2}\int_{\R^2}\langle\nabla g(x),\Gamma\nabla g(x)\rangle dx\bigg\}.\nn 
\end{align} 
We see from (\ref{4.75}) and (\ref{4.76}) that 
\begin{align} \liminf_{t\to\infty}&{1\over b_t}\log 
\E\exp\bigg\{\theta\sqrt{b_t\over t} 
(\log t) 
\vert A_{ t}(\eps )\vert^{1/2}\bigg\}\label{4.77}\\ &\ge\sup_{g\in {\cal 
F}}\bigg\{2\pi\theta\sqrt{\det(\Gamma)} 
\int_{\R^2}f(x)(g^2\ast h_\eps)(x) dx \nn\\ 
&\hspace{ 2in}-{1\over 
2}\int_{\R^2}\langle\nabla g(x),\Gamma\nabla g(x)\rangle dx\bigg\}.\nn 
\end{align} 
  Taking the supremum over $f$ on the right gives 
\begin{align} \liminf_{t\to\infty}&{1\over b_t}\log 
\E\exp\bigg\{\theta\sqrt{b_t\over t} 
(\log t) 
\vert A_{ t}(\eps )\vert^{1/2}\bigg\}\label{4.78}\\ &\ge\sup_{g\in {\cal 
F}}\bigg\{2\pi\theta\sqrt{\det(\Gamma)} 
\bigg(\int_{\R^2}\vert (g^2\ast h_\eps)(x)\vert^2dx\bigg)^{1/2} 
\nn\\ 
&\hspace{ 2in} -{1\over 2}\int_{\R^2}\langle\nabla g(x),\Gamma\nabla g(x)\rangle 
dx\bigg\}.\nn 
\end{align} 
This completes the proof of (\ref{1.00a}). 
 
To prove (\ref{4.145}),  in (\ref{1.00a}) we replace $t$ by 
$2^{-N}t$, 
$\theta$ by 
$2^{-N/2}\theta$, 
$b_t$ by $\wt 
b_t=:b_{ 2^{N}t}$ and $\eps$ by $2^{N/2}\eps$ to find that 
\begin{align} 
\lim_{t\to\infty}&{1\over b_t}\log 
\E\exp\bigg\{\theta\sqrt{b_t\over t} 
(\log t)\label{4.145a}\\ 
&\hskip.5in\times\bigg(\Lambda_\eps \Big({t\over 
b_t}\Big)^{-2}\sum_{x\in \Z^2}\bigg[\sum_{y\in Z[0,2^{-N}t]} 
h_\eps\Big(\sqrt{b_t\over t}(x-y)\Big)\bigg]^2 
\bigg)^{1/2}\bigg\}\nn\\ 
&= 
\lim_{t\to\infty}{1\over \wt b_{2^{-N}t}}\log 
\E\exp\bigg\{2^{-N/2}\theta\sqrt{\wt b_{2^{-N}t}\over 2^{-N}t} 
(\log t)\nn\\ 
&\hskip.4in\times\bigg(\Lambda_{2^{N/2}\eps} \Big({2^{-N}t\over 
\wt b_{2^{-N}t}}\Big)^{-2}\sum_{x\in \Z^2}\bigg[\sum_{y\in Z[0,2^{-N}t]} 
h_{2^{N/2}\eps} \Big(\sqrt{\wt b_{2^{-N}t}\over 
2^{-N}t}(x-y)\Big)\bigg]^2 
\bigg)^{1/2}\bigg\}\nn\\ &=\sup_{g\in {\cal F}}\bigg\{2\pi 
2^{-N/2}\theta\sqrt{\det(\Gamma)} 
\bigg(\int_{\R^2}\vert (g^2\ast h_{2^{N/2}\eps} )(x)\vert^2dx\bigg)^{1/2} 
\nn\\&~~~-{1\over 2}\int_{\R^2}\langle\nabla g(x),\Gamma\nabla g(x)\rangle^2 
dx\bigg\}\nn\\ &\le\sup_{g\in {\cal F}}\bigg\{2\pi 
2^{-N/2}\theta\sqrt{\det(\Gamma)} 
\bigg(\int_{\R^2}\vert g(x)\vert^4dx\bigg)^{1/2} \nn\\ 
&~~~~-{1\over 
2}\int_{\R^2}\langle\nabla g(x),\Gamma\nabla g(x)\rangle^2 dx\bigg\}.\nn\\ 
&=\big(2\pi 
2^{-N/2}\theta\big)^2\sqrt{\det(\Gamma)} 
\sup_{f\in {\cal F}}\bigg\{ 
\bigg(\int_{\R^2}\vert f(x)\vert^4dx\bigg)^{1/2} \nn\\ 
&~~~~-{1\over 
2}\int_{\R^2}\vert\nabla f(x)\vert^2 dx\bigg\}.\nn\\ 
&=2^{-N+2}\pi^2\theta^2\sqrt{\det(\Gamma)}\kappa(2,2)^4,\nn 
\end{align} 
where the third step follows from Jensen's inequality, the fourth step 
follows from the substitution 
$g(x)=\sqrt{\vert\det(A)\vert}f(Ax)$ with 
the $2\times 2$ matrix $A$ satisfying 
$$ 
A^{\tau}\Gamma A=\big(2\pi 
2^{-N/2}\theta\big)^2\sqrt{\det(\Gamma)}I_{2\times 2}, 
$$ 
and the last step follows from Lemma 7.2 in [7]. 
\qed

\section{Exponential approximation}\label{sec-expapp} 
 
Let $t_1,\cdots, t_a\ge 0$ and write 
\begin{equation} 
\label{4.26} 
\Delta_1=[0, t_1],\hskip.1in\hbox{and}\hskip.1in 
\Delta_k=\bigg[\sum_{j=1}^{k-1}t_j,\hskip.1in \sum_{j=1}^kt_j\bigg] 
\hskip.2in k=2,\cdots, a. 
\end{equation} 
Let $p(x)$ be a positive symmetric function on $\Z^2$ with $  \sum_{x\in\Z^2}p(x)=1 $ and define 
\begin{equation} 
\label{4.27} L=\sum_{j,k=1}^a\bigg[\big\vert Z(\Delta_j)\cap 
Z'(\Delta_k)\big\vert -\sum_{x\in \Z^2}p(x)\big\vert Z(\Delta_j)\cap 
\big(Z'(\Delta_k)+x\big) 
\big\vert\bigg], 
\end{equation} 
and 
\begin{align} 
L_j&=\big\vert Z[0,t_j]\cap Z'[0,t_j]\big\vert  
\label{4.28}  
\\&\quad -\sum_{x\in 
\Z^2}p(x)\big\vert Z[0,t_j]\cap \big(Z'[0,t_j]+x\big) 
\big\vert,\hskip.2in j=1,\cdots, a. 
\nn 
\end{align}

\begin{lemma}\label{LU4} For any $m\ge 1$, 
\begin{equation} 
\E L^m\geq 0\label{4.29s} 
\end{equation} 
and 
\begin{equation} 
\label{4.29} 
\Big\{\E L^m\Big\}^{1/2} 
\le\sum_{\scriptstyle k_1 +\cdots +k_a=m 
\atop\scriptstyle k_1,\cdots, k_a\ge 0} {m!\over k_1!\cdots k_a!}\Big\{\E 
\vert L_1\vert^{k_1}\Big\}^{1/2}\cdots 
\Big\{\E \vert L_a\vert^{k_a}\Big\}^{1/2}. 
\end{equation} Consequently, for any $\theta >0$ 
\begin{equation} 
\label{4.30} 
\sum_{m=0}^\infty {\theta^m\over m!}\Big\{\E L^m\Big\}^{1/2} 
\le\prod_{j=1}^a\sum_{m=0}^\infty {\theta^m\over m!} 
\Big\{\E \vert L_j\vert^m\Big\}^{1/2}. 
\end{equation} 
\end{lemma}

\proof Write 
\begin{equation} 
\label{4.31} 
\widehat{p}(\lambda)=\sum_{x\in\Z^2}p(x)e^{i\lambda\cdot x}. 
\end{equation} 
We note that 
\begin{equation} 
|\widehat{p}(\lambda)|\leq \widehat{p}(0)=1.\label{4.31a} 
\end{equation} 
  Notice also that 
\begin{equation} 
\label{4.32} L={1\over (2\pi)^2}\int_{[-\pi,\pi]^2} 
\big[1-\widehat{p}(\lambda)\big]\Big[\sum_{j=1}^a\sum_{x\in 
Z(\Delta_j)} e^{i\lambda\cdot x}\Big]\Big[\sum_{j'=1}^a\,\sum_{x'\in 
Z'(\Delta_{j'})} e^{-i\lambda\cdot x'}\Big]\,d\lambda . 
\end{equation} 
  We therefore have 
\begin{align} 
\E L^m&={1\over (2\pi)^{2m}}\int_{([-\pi,\pi]^2)^m} 
\bigg\vert \E\prod_{k=1}^m\sum_{j=1}^a\sum_{x_{ k}\in Z(\Delta_j)} 
e^{i\lambda_{ k}\cdot x_{ k}}\bigg\vert^2 
\Big(\prod_{k=1}^m\big[1-\widehat{p}(\lambda_k)\big]d\lambda_k\Big)\nn\\ 
&={1\over (2\pi)^{2m}}\int_{([-\pi,\pi]^2)^m}\bigg\vert  
\sum_{l_1,\cdots,l_m=1}^a 
\E\Big(H_{l_1}(\lambda_1)\cdots 
H_{l_m}(\lambda_m)\Big)\bigg\vert^2 
\nn\\ 
&\hspace{3in} 
\Big(\prod_{k=1}^m\big[1-\widehat{p}(\lambda_k)\big]d\lambda_k\Big),\label{4.33} 
\end{align} where 
\begin{equation} 
\label{4.34} H_j(\lambda)=\sum_{x\in Z(\Delta_j)}e^{i\lambda\cdot x}. 
\end{equation} 
This proves (\ref{4.29s}) and implies that 
\begin{align} 
\Big\{\E& L^m\Big\}^{1/2} 
\label{4.34A}\\ 
&\le {1\over 
(2\pi)^m}\sum_{l_1,\cdots,l_m=1}^a 
\bigg\{\int_{([-\pi,\pi]^2)^m} 
\Big\vert \E\Big(H_{l_1}(\lambda_1)\cdots 
H_{l_m}(\lambda_m)\Big)\Big\vert^2 
\nn \\ 
&\hspace{2in} \Big(\prod_{k=1}^m\big[1-\widehat{p} 
(\lambda_k)\big]d\lambda_k\Big)\bigg\}^{1/2}.\nn 
\end{align} 
Note that for any $k>j$ we can write 
\begin{equation} 
\label{4.34a} H_k(\lambda)=\sum_{x\in Z(\Delta_k)}e^{i\lambda\cdot x} 
=e^{i\lambda\cdot Z(t_j)}H^{ ( j)}_k(\lambda), 
\end{equation} 
where 
\begin{equation} 
\label{4.34b} H^{ ( j)}_k(\lambda)=\sum_{x\in 
Z(\Delta_k)-Z(t_j)}e^{i\lambda\cdot x} 
\end{equation} 
is independent of $\mathcal{F}_{ t_j}$. 
 
Let $l_1,\cdots, l_m$ be fixed and let $k_j=\sum_{ i=1}^{ m}\de(l_i,j 
)$  be  the number of $l$'s which are equal to $j$, for each $1\leq j\leq  a$. 
Then using independence 
\begin{align}  
&\int_{([-\pi,\pi]^2)^m} 
\Big\vert \E\Big(H_{l_1}(\lambda_1)\cdots 
H_{l_m}(\lambda_m)\Big)\Big\vert^2 
\Big(\prod_{k=1}^m\big[1-\widehat{p}(\lambda_k)\big]d\lambda_k\Big) 
\label{4.35}\\ &=\int_{([-\pi,\pi]^2)^m} 
\Big\vert \E\prod_{j=1}^a\Big(H_j(\lambda_{j,1})\cdots 
H_j(\lambda_{j,k_j})\Big)\Big\vert^2 
\Big(\prod_{j=1}^a\prod_{l=1}^{k_j} 
\big[1-\widehat{p}(\lambda_{j,l})\big]d\lambda_{j,l}\Big) 
\nn\\ &=\int_{([-\pi,\pi]^2)^m} 
\Big\vert \E\bigg[\exp\Big\{i\Big(\sum_{j=2}^a\sum_{l=1}^{k_j} 
\lambda_{j,l}\Big)\cdot Z(t_1)\Big\}\times 
\Big(H_1(\lambda_{1,1})\cdots  H_1(\lambda_{1,k_1})\Big)\bigg]\nn\\ & 
\hspace{ 1in}\E\Big(\prod_{j=2}^a\Big(H^{ ( 1)}_j(\lambda_{j,1})\cdots 
H^{ ( 1)}_j(\lambda_{j,k_j})\Big)\Big\vert^2 
\Big(\prod_{j=1}^a\prod_{l=1}^{k_j} 
\big[1-\widehat{p}(\lambda_{j,l})\big]d\lambda_{j,l}\Big)\nn\\ 
&=\int_{([-\pi,\pi]^2)^{m-k_1}} 
\Big\vert \E\Big(\prod_{j=2}^a\Big(H^{ ( 1)}_j(\lambda_{j,1})\cdots 
H^{ ( 1)}_j(\lambda_{j,k_j})\Big)\Big\vert^2\nn\\ 
&~~~~F(\lambda_{2,1},\cdots,\lambda_{2,k_2};\cdots;\lambda_{a,1}, \cdots,\lambda_{a,k_a}) 
\Big(\prod_{j=2}^a\prod_{l=1}^{k_j} 
\big[1-\widehat{p}(\lambda_{j,l})\big]d\lambda_{j,l}\Big),\nn            
\end{align} 
where  
\begin{align}  
&F(\lambda_{2,1},\cdots,\lambda_{2,k_2};\cdots;\lambda_{a,1}, \cdots,\lambda_{a,k_a}) 
\label{4.35a}\\ &=\int_{([-\pi,\pi]^2)^{k_1}} 
\Bigg\vert\E\bigg[\exp\Big\{i\Big(\sum_{j=2}^a\sum_{l=1}^{k_j} 
\lambda_{j,l}\Big)\cdot Z(t_1)\Big\}\times\Big(H_1(\lambda_{1,1})\cdots 
H_1(\lambda_{1,k_1})\Big)\bigg]\Bigg\vert^2 \nonumber\\ 
&\hspace{3in} 
\Big(\prod_{l=1}^{k_1} 
\big[1-\widehat{p}(\lambda_{1,l})\big]d\lambda_{1,l}\Big).\nn            
\end{align} 
Notice that by symmetry 
\begin{equation} 
\label{4.36} 
\E\bigg[\exp\Big\{i\Big(\sum_{j=2}^a\sum_{l=1}^{k_j} 
\lambda_{j,l}\Big)\cdot Z(t_1)\Big\}\Big(H_1(\lambda_{1,1})\cdots 
H_1(\lambda_{1,k_1})\Big)\bigg] 
\end{equation} 
  is real valued. Hence if $Z'$ denotes an independent copy of $Z$, and $H'_{ 
1}$ is obtained from $H_{ 1}$ by replacing  $Z$ by  $Z'$, 
\begin{eqnarray} &&\qquad  
F(\lambda_{2,1},\cdots,\lambda_{2,k_2};\cdots;\lambda_{a,1}, \cdots,\lambda_{a,k_a}) \label{4.37}\\  
&&=\int_{([-\pi,\pi]^2)^{k_1}} 
\E\bigg[\exp\Big\{i\Big(\sum_{j=2}^a\sum_{l=1}^{k_j} 
\lambda_{j,l}\Big)\cdot 
\big(Z(t_1)+Z'(t_1)\big)\Big\}\nn\\ 
&&\hspace{.5in}\times\prod_{l=1}^{k_1}\Big(H_1(\lambda_{1,l})H_1'( 
\lambda_{1,l}) 
\Big)\bigg] 
\Big(\prod_{l=1}^{k_1} 
\big[1-\widehat{p}(\lambda_{1,l})\big]d\lambda_{1,l}\Big)\nn\\ 
&&=\E\bigg[\exp\Big\{i\Big(\sum_{j=2}^a\sum_{l=1}^{k_j} 
\lambda_{j,l}\Big)\cdot 
\big(Z(t_1)+Z'(t_1)\big)\Big\}\nn\\ 
&&\hspace{ .4in}\times\int_{([-\pi,\pi]^2)^{k_1}} 
\prod_{l=1}^{k_1}\Big(H_1(\lambda_{1,l})H_1'( \lambda_{1,l}) 
\Big)\bigg] 
\Big(\prod_{l=1}^{k_1} 
\big[1-\widehat{p}(\lambda_{1,l})\big]d\lambda_{1,l}\Big).\nn 
\end{eqnarray} 
By the fact that 
\begin{eqnarray}&&\int_{([-\pi,\pi]^2)^{k_1}} 
\prod_{l=1}^{k_1}\Big(H_1(\lambda_{1,l})H_1'( \lambda_{1,l}) 
\Big) 
\Big(\prod_{l=1}^{k_1} 
\big[1-\widehat{p}(\lambda_{1,l})\big]d\lambda_{1,l}\Big)\label{4.38}\\ 
&&=\bigg[\int_{[-\pi,\pi]^2}\big[1-\widehat{p}(\lambda)\big] 
H_1(\lambda)H_1'(\lambda)\,d\lambda\bigg]^{k_1} =(2\pi)^{2k_1} 
L_1^{k_1},\nn 
\end{eqnarray} 
we have proved that 
\begin{eqnarray} &&\int_{([-\pi,\pi]^2)^m} 
\Big\vert \E\Big(H_{l_1}(\lambda_1)\cdots 
H_{l_m}(\lambda_m)\Big)\Big\vert^2 
\Big(\prod_{k=1}^m\big[1-\widehat{p}(\lambda_k)\big]d\lambda_k\Big) 
\nn\\ &&\le (2\pi)^{2k_1} 
\E\vert L_1\vert^{k_1} 
\int_{([-\pi,\pi]^2)^{m-k_1}}\Big\vert 
\E\Big(\prod_{j=2}^a\Big(H^{ ( 1)}_j(\lambda_{j,1})\cdots 
H^{ ( 1)}_j(\lambda_{j,k_j})\Big)\Big\vert^2\nn\\ &&\hskip .3in 
\Big(\prod_{j=2}^a\prod_{l=1}^{k_j} 
\big[1-\widehat{p}(\lambda_{j,l})\big]d\lambda_{j,l}\Big).\label{4.39} 
\end{eqnarray} 
  Repeating the above procedure, 
\begin{eqnarray}&&\int_{([-\pi,\pi]^2)^m} 
\Big\vert \E\Big(H_{l_1}(\lambda_1)\cdots 
H_{l_m}(\lambda_m)\Big)\Big\vert^2 
\Big(\prod_{k=1}^m\big[1-\widehat{p}(\lambda_k)\big]d\lambda_k\Big) 
\nn\\ 
&&\hspace{ 1in}\le\prod_{j=1}^a\Big\{(2\pi)^{2k_j}\E\vert 
L_j\vert^{k_j}\Big\} =(2\pi)^{2m}\prod_{j=1}^a\E\vert L_j\vert^{k_j}. 
\label{4.40} 
\end{eqnarray}
Our Lemma now follows from (\ref{4.34A}). 
\qed

{\bf  Proof of Lemma \ref{LU7}}:  Define 
\begin{equation} q_{t,\eps}(x)=\Lambda_\eps \Big({t\over 
b_t}\Big)^{-2}\sum_{z\in\Z^2} h_\eps\Big(\sqrt{b_t\over t}(x-z)\Big) 
h_\eps\Big(\sqrt{b_t\over t}z\Big),\hskip.2in x\in\Z^2. 
\label{4.82} 
\end{equation} 
Then $q_{t,\eps}(x)$ is a probability density on $\Z^2$. We claim that 
\begin{equation} 
B^{ ( 0)}_{ t}( \eps)=\sum_{x\in\Z^2}q_{t,\eps}(x) 
\big\vert Z[0, t]\cap \big(x+Z'[0, t]\big)\big\vert.\label{4.82r} 
\end{equation} 
This follows from the fact  that 
\begin{eqnarray} &&\sum_{x\in\Z^2}\sum_{y\in Z[0,t]} 
h_\eps\Big(\sqrt{b_t\over t}(x-y)\Big) 
\sum_{y'\in Z'[0,t]} h_\eps\Big(\sqrt{b_t\over 
t}(x-y')\Big)\label{4.68a}\\ && 
=\sum_{x\in\Z^2}\sum_{y'\in Z'[0,t]} h_\eps\Big(\sqrt{b_t\over 
t}x\Big) 
\sum_{y\in \Z^2} h_\eps\Big(\sqrt{b_t\over 
t}(x+y'-y)\Big)1_{\{y\in Z[0,t]\}}\nonumber\\ && 
=\sum_{x\in\Z^2}\sum_{y'\in Z'[0,t]} h_\eps\Big(\sqrt{b_t\over 
t}x\Big) 
\sum_{y\in\Z^2} h_\eps\Big(\sqrt{b_t\over 
t}(x-y)\Big)1_{\{y+y'\in Z[0,t]\}}\nonumber\\ && 
=\sum_{y\in\Z^2}\sum_{x\in\Z^2} h_\eps\Big(\sqrt{b_t\over 
t}x\Big) h_\eps\Big(\sqrt{b_t\over 
t}(x-y)\Big) 
\big\vert Z'[0,t]\cap \big(Z[0, 
t]-y\big)\big\vert\nonumber 
\end{eqnarray} 
and 
\begin{equation} 
\big\vert Z'[0,t]\cap \big(Z[0, 
t]-y\big)\big\vert=\big\vert  Z[0,t]\cap\big(Z'[0,t]+y\big) 
\big\vert.\label{4.68b} 
\end{equation}

  Write $\gamma_t=t/[b_t]$ and 
$\Delta_j=[(j-1)\gamma_t,j\gamma_t]$, $j=1,\cdots, [b_t]$. 
Note that 
\begin{eqnarray} &\sum_{j=1}^{[b_t]}\big\vert Z(\Delta_j)\cap 
Z'[0,t]\big\vert -\sum_{1\le j< k\le [b_t]}\big\vert Z(\Delta_j)\cap 
Z(\Delta_k)\cap Z'[0,t]\big\vert\nn\\ &\le\big\vert Z[0, t]\cap Z'[0, 
t]\big\vert 
\le\sum_{j=1}^{[b_t]}\big\vert Z(\Delta_j)\cap Z'[0,t]\big\vert\label{1.11} 
\end{eqnarray} 
and similarly 
\begin{align} &\sum_{j=1}^{[b_t]}\sum_{x\in\Z^2}q_{t,\eps}(x) 
\big\vert Z(\Delta_j)\cap \big(x+Z'[0,t]\big)\big\vert\label{1.12}\\ 
&-\sum_{1\le j< k\le [b_t]}\sum_{x\in\Z^2}q_{t,\eps}(x) 
\big\vert Z(\Delta_j)\cap Z(\Delta_k)\cap \big(x+Z'[0,t]\big)\big\vert\nn\\ 
&\le \sum_{x\in\Z^2}q_{t,\eps}(x) 
\big\vert Z[0, t]\cap \big(x+Z'[0, t]\big)\big\vert\nn\\ &\le 
\sum_{j=1}^{[b_t]}\sum_{x\in\Z^2}q_{t,\eps}(x) 
\big\vert Z(\Delta_j)\cap \big(x+Z'[0,t]\big)\big\vert.\nn 
\end{align} Hence, 
\begin{align}&\bigg\vert \big\vert Z[0, t]\cap Z'[0, 
t]\big\vert-\sum_{x\in\Z^2}q_{t,\eps}(x) 
\big\vert Z[0, t]\cap \big(x+Z'[0, t]\big)\big\vert \bigg\vert\nn\\ 
& 
\le \bigg\vert \sum_{j=1}^{[b_t]}\Big[\big\vert Z(\Delta_j)\cap 
Z'[0,t]\big\vert- 
\sum_{x\in\Z^2}q_{t,\eps}(x) 
\big\vert Z(\Delta_j)\cap \big(x+Z'[0,t]\big)\big\vert\Big]\bigg\vert\nn\\ 
&+\sum_{1\le j< k\le [b_t]}\big\vert Z(\Delta_j)\cap Z(\Delta_k)\cap 
Z'[0,t]\big\vert\nn\\ &+\sum_{1\le j< k\le 
[b_t]}\sum_{x\in\Z^2}q_{t,\eps}(x) 
\big\vert Z(\Delta_j)\cap Z(\Delta_k)\cap 
\big(x+Z'[0,t]\big)\big\vert.\label{1.13} 
\end{align} 
 
We first take care of the last two terms. This is the easy step. Write 
\begin{align} &\eta(t,\eps)=\sum_{1\le j< k\le [b_t]}\big\vert 
Z(\Delta_j)\cap Z(\Delta_k)\cap Z'[0,t]\big\vert\label{1.14}\\ 
&+\sum_{1\le j< k\le [b_t]}\sum_{x\in\Z^2}q_{t,\eps}(x) 
\big\vert Z(\Delta_j)\cap Z(\Delta_k)\cap \big(x+Z'[0,t]\big)\big\vert.\nn 
\end{align}

It follows from (\ref{4.24}) that 
\begin{equation} 
\sup_{ t,j,k,x}\E\exp\bigg\{c{(\log t)^{3/2}\over \sqrt{t}} 
\big\vert Z(\Delta_j)\cap Z(\Delta_k)\cap 
\big(x+Z'[0,t]\big)\big\vert^{1/2}\bigg\}<\infty. 
\end{equation} 
for some $c>0$. 
Hence, if $b_t=o\big((\log 
t)^{1/5}\big)$, then for any $\theta >0$ we can find $t_{ 0}<\ff$ such that 
\begin{align} \sup_{ t\geq t_{ 0}}&\E\exp\bigg\{\theta\sqrt{b_t\over 
t} (\log t)\eta(t,\eps)^{1/2}\bigg\}\label{1.15}\\ &\le\sup_{ t\geq t_{ 0}}\, 
\sup_{ j,k,x}\E\exp\bigg\{\theta\sqrt{b_t\over t}(\log t)b_t^2 
\big\vert Z(\Delta_j)\cap Z(\Delta_k)\cap 
\big(x+Z'[0,t]\big)\big\vert^{1/2}\bigg\}\nn\\ 
&<\infty.\nn            
\end{align} 
  Hence 
\begin{equation} 
\limsup_{t\to\infty}{1\over b_t}\log 
\E\exp\Big\{\theta\sqrt{b_t\over t}(\log t) 
\eta(t,\eps)^{1/2}\Big\}=0. 
\label{4.86} 
\end{equation} 
 
\medskip To handle the first term on the right hand side of (\ref{1.13}) set 
\begin{equation} 
\xi(t,\eps)=\sum_{j=1}^{[b_t]}\Big[\big\vert Z(\Delta_j)\cap 
Z'[0,t]\big\vert- 
\sum_{x\in\Z^2}q_{t,\eps}(x) 
\big\vert Z(\Delta_j)\cap \big(x+Z'[0,t]\big)\big\vert\Big]. 
\label{4.87} 
\end{equation} 
Using Fubini, independence and then the Cauchy-Schwarz inequality we have 
\begin{align} 
\big\vert&\E\xi^m (t,\eps)\big\vert=\label{1.16}\\ 
&(2\pi)^{-2m} 
\bigg\vert\E\int_{([-\pi,\pi]^2)^m} 
\Big(\prod_{k=1}^m 
\big[1-\widehat{q}_{t,\eps}(\lambda_{ k})\big]\Big) 
\nn\\ &\hspace{ 
.2in}\times\bigg[\prod_{k=1}^m\sum_{x'_{ k}\in Z'[0,t]}e^{i\lambda_{ 
k}\cdot x'_{ k}}\bigg] 
\bigg[\prod_{k=1}^m\sum_{j=1}^{[b_n]}\sum_{x_{ k}\in Z(\Delta_j)} 
e^{-i\lambda_{ k}\cdot x_{ k}} 
\bigg]\,d\lambda _1\cdots d\lambda_m\bigg\vert\nn\\ 
&\le(2\pi)^{-2m}\bigg\{\int_{([-\pi,\pi]^2)^m} 
\Big(\prod_{k=1}^m \big[1-\widehat{q}_{t,\eps}(\lambda_{ k})\big] 
\Big)\bigg\vert\E\prod_{k=1}^m\sum_{x_{ k}\in Z[0,t]} e^{i\lambda_{ k}\cdot 
x_{ k}}\bigg\vert^2\,d\lambda _1\cdots d\lambda_m\bigg\}^{1/2}\nn\\ 
&\times\bigg\{\int_{([-\pi,\pi]^2)^m} 
\Big(\prod_{k=1}^m \big[1-\widehat{q}_{t,\eps}(\lambda_{ k})\big] 
\Big)\bigg\vert\E\prod_{k=1}^m\sum_{j=1}^{[b_t]}\sum_{x_{ k}\in 
Z(\Delta_j)} e^{i\lambda_{ k}\cdot x_{ k}}\bigg\vert^2\,d\lambda _1\cdots 
d\lambda_m\bigg\}^{1/2}\nn\\ &\le\Big\{\E \vert Z[0,t]\cap 
Z'[0,t]\vert^m\Big\}^{1/2} 
\Big\{\E\zeta^m(t,\eps)\Big\}^{1/2},\nn 
\end{align} 
  where 
\begin{equation} 
\zeta (t,\eps)=\sum_{j,k =1}^{[b_t]}\Big[\big\vert Z(\Delta_j)\cap 
Z'(\Delta_k)\big\vert- 
\sum_{x\in\Z^2}q_{t,\eps}(x) 
\big\vert Z(\Delta_j)\cap \big(x+Z'(\Delta_k)\big)\big\vert\Big] 
\label{4.88} 
\end{equation} 
and we have used the fact that 
$1-\widehat{q}_{t,\eps}(\lambda)\le 1$ in the last step. 
Note that in the notation of (\ref{4.27}), $\zeta (t,\eps)=L$  with $p( 
x)=q_{t,\eps}( x)$, so that by (\ref{4.29s}), for all $m\geq 1$ 
\begin{equation} 
\E\zeta^m(t,\eps)\ge 0.\label{4.29ss} 
\end{equation} 
 
Let $\delta >0$ be fixed for a while. By Cauchy-Schwarz and then (\ref{1.16}) 
\begin{eqnarray} && 
\E\cosh\Big\{\theta \sqrt{b_t\over t}(\log t)\vert\xi(t,\eps 
)\vert^{1/2}\Big\}\label{4.99}\\ && 
=\sum_{m=0}^\infty{\theta^{2m}\over (2m)!}\Big(\sqrt{b_t\over t} 
(\log t)\Big)^{2m}\E|\xi^m(t,\eps)|\nonumber\\ 
&&\leq 
\sum_{m=0}^\infty{\theta^{2m}\over (2m)!}\Big(\sqrt{b_t\over t} (\log 
t)\Big)^{2m} 
\Big\{\E\xi^{2m}(t,\eps)\Big\}^{1/2}\nn\\ 
&&\le\bigg\{\sum_{m=0}^\infty{(\delta\theta)^{2m}\over (2m)!} 
\Big(\sqrt{b_t\over t}(\log t)\Big)^{2m} 
\Big\{\E\vert Z[0,t]\cap Z'[0,t]\vert^{2m}\Big\}^{1/2}\bigg\}^{1/2}\nn\\ 
&&\times 
\bigg\{\sum_{m=0}^\infty{(\delta^{-1}\theta)^{2m}\over (2m)!} 
\Big(\sqrt{b_t\over t}(\log t)\Big)^{2m} 
\Big\{\E\zeta^{2m}(t,\eps)\Big\}^{1/2}\bigg\}^{1/2}\nn\\ 
&&\le\bigg\{\sum_{m=0}^\infty{(\delta\theta)^m\over m!} 
\Big(\sqrt{b_t\over t}(\log t)\Big)^m 
\Big\{\E\vert Z[0,t]\cap Z'[0,t]\vert^m\Big\}^{1/2}\bigg\}^{1/2}\nn\\ 
&&\times 
\bigg\{\sum_{m=0}^\infty{(\delta^{-1}\theta)^m\over m!}\Big(\sqrt{b_t\over 
t}(\log t)\Big)^m 
\Big\{\E\zeta^m(t,\eps)\Big\}^{1/2}\bigg\}^{1/2},\nn 
\end{eqnarray} 
  where in the last step we used  (\ref{4.29ss}) and the fact that 
$\vert Z[0,t]\cap Z'[0,t]\vert\ge 0.$ 
 
By \cite[(2.11)]{Cb}, there is a $C>0$ independent of $\delta$ and $\theta$ 
such that 
\begin{equation} 
\lim_{t\to\infty}{1\over b_t}\log\sum_{m=0}^\infty{(\delta \theta)^m\over m!} 
\Big(\sqrt{b_t\over t}(\log t)\Big)^m 
\Big\{\E\vert Z[0,t]\cap Z'[0,t]\vert^m\Big\}^{1/2} 
= C(\delta\theta)^2. 
\label{4.90} 
\end{equation} 
In addition, by Lemma \ref{LU4} 
\begin{eqnarray} & 
\sum_{m=0}^\infty{(\delta^{-1}\theta)^m\over m!}\Big(\sqrt{b_t\over t}(\log 
t)\Big)^m 
\Big\{\E\zeta^m(t,\eps)\Big\}^{1/2}\label{1.18}\\ &\le 
\bigg\{\sum_{m=0}^\infty{(\delta^{-1}\theta)^m\over m!}\Big(\sqrt{b_t\over 
t}(\log t)\Big)^m 
\Big\{\E \vert\beta(t,\eps)\vert^m\Big\}^{1/2}\bigg\}^{[b_t]},\nn 
\end{eqnarray} 
  where 
\begin{equation} 
\beta (t,\eps)=\big\vert Z[0,\gamma_t]\cap Z'[0,\gamma_t]\big\vert 
-\sum_{x\in\Z^2}q_{t,\eps}(x)\big\vert Z[0,\gamma_t]\cap 
\big(x+Z'[0,\gamma_t])\big\vert. 
\label{4.91} 
\end{equation} 
Recall that $q_{t,\eps}(x)$ is defined by (\ref{4.82}) and 
$\gamma_t=t/[b_t]$. As in the proof of (\ref{4.82r}) we can check that
\[\sum_{x\in\Z^2}q_{t,\eps}(x)\big\vert Z[0,\gamma_t]\cap 
\big(x+Z'[0,\gamma_t])\big\vert=B_{ \gamma_t,1},\] 
see (\ref{5.01b}). By 
Lemma 
\ref{LU5} (with $t$ replaced by 
$\gamma_t$), 
\begin{equation} {b_t(\log t)^2\over t}\beta (t,\eps)\buildrel 
d\over\longrightarrow (2\pi)^2\det(\Gamma)\Big[\alpha 
([0,1]^2)-\alpha_\eps ([0,1]^2)\Big]. 
\label{4.92} 
\end{equation} 
  By Lemma \ref{LU2} (with $p=2$), 
\begin{equation} 
\E \vert\beta(t,\eps)\vert^m\le 2\sup_{ x}\E^{( 0,x)}\vert Z[0,\gamma_t]\cap Z'[0,\gamma_t]\vert^m 
\le m!C^m \Big({t\over b_t}(\log t)^{-2}\Big)^m. 
\label{4.93} 
\end{equation} 
Hence, 
\begin{eqnarray} &\lim_{t\to\infty} 
\sum_{m=0}^\infty{(\delta^{-1}\theta)^m\over m!}\Big(\sqrt{b_t\over t}(\log 
t)\Big)^m 
\Big\{\E \vert\beta(t,\eps)\vert^m\Big\}^{1/2}\label{1.19}\\ 
&=\sum_{m=0}^\infty{(\delta^{-1}\theta)^m\over 
m!}\Big((2\pi)\sqrt{\det(\Gamma)}\Big)^m 
\Big\{\E\Big\vert\alpha ([0,1]^2) -\alpha_\eps 
([0,1]^2)\Big\vert^m\Big\}^{1/2}.\nn 
\end{eqnarray} 
So by (\ref{1.18}) we have 
\begin{eqnarray} &\limsup_{t\to\infty}{1\over b_t}\log 
\sum_{m=0}^\infty{(\delta^{-1}\theta)^m\over m!}\Big(\sqrt{b_t\over t}(\log 
t)\Big)^m 
\Big\{\E\zeta^m(t,\eps)\Big\}^{1/2}\label{1.20}\\ &\le 
\log \sum_{m=0}^\infty{(\delta^{-1}\theta)^m\over 
m!}\Big((2\pi)\sqrt{\det(\Gamma)}\Big)^m 
\Big\{\E\Big\vert\alpha ([0,1]^2) -\alpha_\eps 
([0,1]^2)\Big\vert^m\Big\}^{1/2}.\nn 
\end{eqnarray} 
  By \cite[Theorem 1, p.183]{LeGall-St.Flour}, 
\begin{equation} 
\E\Big\vert\alpha ([0,1]^2) -\alpha_\eps 
([0,1]^2)\Big\vert^m\longrightarrow 0\hskip.2in \mbox{as } 
\eps\to 0^+, 
\label{4.94} 
\end{equation} 
  for all $m\ge 1$. In addition, by \cite[(1.12)]{C},  there is a constant 
$C>0$ such that 
\begin{equation} 
\E\Big\vert\alpha ([0,1]^2) -\alpha_\eps 
([0,1]^2)\Big\vert^m\le\E\alpha^m ([0,1]^2)\le m! C^m 
\label{4.95} 
\end{equation} 
   for all $m\ge 1$. By dominated convergence, therefore, 
\begin{equation} 
\sum_{m=0}^\infty{(\delta^{-1}\theta)^m\over 
m!}\Big((2\pi)\sqrt{\det(\Gamma)}\Big)^m 
\Big\{\E\Big\vert\alpha ([0,1]^2) -\alpha_\eps 
([0,1]^2)\Big\vert^m\Big\}^{1/2} 
\longrightarrow 1 
\label{4.96} 
\end{equation} 
as $\eps\to 0^+$. 
(Alternatively, this follows immediately from \cite[(6.29)]{CR}). 
  Thus 
\begin{equation} 
\lim_{\eps\to 0^+}\limsup_{t\to\infty}{1\over b_t}\log 
\sum_{m=0}^\infty{(\delta^{-1}\theta)^m\over m!}\Big(\sqrt{b_t\over t}(\log 
t)\Big)^m 
\Big\{\E\zeta^m(t,\eps)\Big\}^{1/2}=0. 
\label{4.97} 
\end{equation}

Summarizing what we have, 
\begin{equation} 
\limsup_{\eps\to 0^+}\limsup_{t\to\infty}{1\over b_t}\log 
\E\cosh\Big\{\theta \sqrt{b_t\over t}(\log t)\vert\xi(t,\eps 
)\vert^{1/2}\Big\}\le C(\delta\theta)^2. 
\label{4.98A} 
\end{equation} 
Letting $\delta\to 0^+$ gives  
\begin{equation} 
\limsup_{\eps\to 0^+}\limsup_{t\to\infty}{1\over b_t}\log 
\E\cosh\Big\{\theta \sqrt{b_t\over t}(\log t)\vert\xi(t,\eps 
)\vert^{1/2}\Big\}=0. 
\label{4.98} 
\end{equation} 
Since $\exp ( x)\leq 2\cosh( x)$ we see from (\ref{4.98}) that 
\begin{equation} 
\limsup_{\eps\to 0^+}\limsup_{t\to\infty}{1\over b_t}\log 
\E\exp\Big\{\theta \sqrt{b_t\over t}(\log t)\vert\xi(t,\eps 
)\vert^{1/2}\Big\} =0. 
\label{4.100} 
\end{equation}

By  (\ref{1.13}) and (\ref{4.82r}) we have thus 
completed the proof of Lemma \ref{LU7} when $j=0$. 
If in (\ref{4.81}) with $j=0$ we replace $t$ by $2^{-j}t$, $\theta$ by 
$2^{-j/2}\theta$, 
$b_t$ by $\wt 
b_t=:b_{ 2^{j}t}$ and $\eps$ by $2^{j/2}\eps$, we obtain (\ref{4.81}) for 
any $j$ (compare the proof of (\ref{4.145})). 
\qed

\section{Laws of the iterated logarithm} 
 
We first prove some lemmas in preparation for the proof
of Theorem \ref{theo-lilrw}. Define
$$\wt \phi_j=\frac{j}{\HH(j)}, \qquad \wt G_j=(R_j-\wt \phi_j)
\frac{\log^2 n}{n},$$
and $K=[\log \log n]+1$.

\begin{lemma}\label{lemsecr21}
There exists a constant $c_1$ such that if $A$ and $B$ are positive
integers and $C=A+B$, then
$$|\wt \phi_C-\wt \phi_A-\wt \phi_B|\leq c_1 \frac{(A\land B)^{1/2}}
{C^{1/2} \log^2 C}.$$
\end{lemma}

\proof The cases when $A$ or $B$ equal 1 are easy, so we suppose $A,B>1$.
Write
$$\wt \phi_C-\wt\phi_A-\wt \phi_B=\frac{C}{\HH(C)} \Big[
-\frac{A}{C} \frac{\HH(C)-\HH(A)}{\HH(A)}
-\frac{B}{C} \frac{\HH(C)-\HH(B)}{\HH(B)}.\Big]$$
By (\ref{secr11}) and (\ref{secr12}), the right hand side is
bounded in absolute value by
$$c_2\frac{C}{\log C} \Big[
-\frac{A}{C} \frac{\log C-\log A}{\log A} 
-\frac{B}{C} \frac{\log C-\log B}{\log B} \Big]
=c_2|\phi_C-\phi_A-\phi_B|,$$
where $\phi_j=j/\log j$. Our result now follows by Lemma 4.2 of \cite{BK}.
\qed

\begin{lemma} \label{lemsecr22}
There exists $\lambda_0$ such that if
$\lambda\geq \lambda_0$, then
$$\P(\max_{m\leq n} \ol R_m>\lambda n \log\log\log n/\log^2 n)\leq (\log n)^{-2}.$$
\end{lemma}

\proof  
Using Lemma \ref{lemsecr21} in place of Lemma 4.2 of \cite{BK} and
with $\wt \phi_j, \wt G_j$ replacing $\phi_j, G_j$, resp., we have by
\cite{BK}, Lemma 4.3 and the proof of Proposition 4.1 (up through the
display in the middle of p. 1390), that
$$\P(\max_{m\leq n} \wt G_m>A\log\log\log n)\leq (\log n)^{-2}$$
if $A$ is large enough. By
(\ref{lrth6.9}) and (\ref{secr11}), we see that
$$\max_{m\leq n} |\ol R_m - (R_m-\wt \phi_m)|\leq c_1\frac{n}{\log^2 n}
=o(n\log\log \log n/\log^2 n),$$
and our result now follows immediately.
\qed

\noindent {\bf Proof of Theorem \ref{theo-lilrw}:} Let $\xi =2\pi 
\sqrt{\det \Gamma}$. We begin with the  upper bound. Let $\eta, 
\eps>0$ be small and let 
$q>1$ be very close to 1. Let $t_i=[q^i]$. If 
$$A_i=\left\{ \ol R_{t_i}\geq (1+\eta) \xi t_i\log\log\log t_i/ 
\log^2 t_i\right\},$$ then it follows from Theorem \ref{theo-range} that 
$\sum_i\P(A_i)<\infty$, and so by Borel-Cantelli, 
$\P(A_i \hbox{ i.o.})=0.$ 
 
\medskip Next, if $\lam$ is sufficiently large, 
\begin{equation} 
\P(\max_{m\leq n} \ol R_m> \lam n\log \log \log n/\log^2 n) 
\leq(\log n)^{-2};\label{erw51} 
\end{equation} 
by Lemma \ref{lemsecr22}.
Let 
$$B_i=\left\{ \max_{t_i\leq k\leq t_{i+1}} [\ol R_k-\ol R_{t_i}] >\eps 
t_i\log 
\log \log t_i/\log ^2 t_i\right\}.$$ By subadditivity $R_k-R_{t_i}\leq 
R_{k-t_i}\circ \theta_{t_i}$, where $\theta_{t_i}$ is the usual shift 
operator of Markov theory. By   Lemma \ref{lemsecr21}
\begin{equation}
\E R_k-\E R_{t_i}\geq \E R_{k-t_i}- c\frac{t_i}{\log ^2 t_i}.\label{6.1a }
\end{equation} 
 So by 
the Markov property, and using the fact that the $\P^x$ law of $R_{k-t_i}$ does not depend on $x$, for $i$ large
\begin{eqnarray}
&& \P(B_i)\label{ 6.1 }\\
&&=\P ( \max_{t_i\leq k\leq t_{i+1}} [R_k- R_{t_i}-(\E R_k-\E R_{t_i} )] >\eps 
t_i\log 
\log \log t_i/\log ^2 t_i)
\nonumber\\
&& \leq  \P ( \max_{t_i\leq k\leq t_{i+1}} [R_k- R_{t_i}-\E R_{k-t_i}] +c\frac{t_i}{\log ^2 t_i}>\eps 
t_i\log 
\log \log t_i/\log ^2 t_i) \nonumber\\
&& \leq  \P^{S_{t_i}} ( \max_{t_i\leq k\leq t_{i+1}} [\ol R_{k-t_i}] >\eps 
t_i\log\log \log t_i/\log ^2 t_i-c\frac{t_i}{\log ^2 t_i}) \nonumber\\
&& \leq \P(\max_{k\leq t_{i+1}-t_i} \ol R_k\geq \tfrac{\eps}{2} t_i 
\log\log\log t_i /\log^2 t_i) \nonumber
\end{eqnarray} 
 If $q$ is sufficiently small, then $\sum_i 
\P(B_i)$ will be summable by (\ref{erw51}). So with probability one, for 
$i$ large enough 
$$\max_{k\leq t_{i+1}} \ol R_k \leq ((1+\eta)\xi+\eps)q t_i \log 
\log \log t_i/\log^2 t_i.$$ Since $\eta$ and $\eps$ are arbitrary, and we 
can take $q$ as close to 1 as we like, this implies the upper bound.

Let $\eta>0$, $t_i=[\exp(i^{1+\frac{\eta}{2}})]$, $V_i=\# S((t_i,t_{i+1}])$,     and set $$C_i= \left\{ \ol V_i > (1-\eta)\xi (t_{i+1}-t_i)\log\log\log (t_{i+1}-t_i)/\log^2 (t_{i+1}-t_i)\right\}.$$ Note that the events $C_i$ are independent. By Theorem \ref{theo-range} and Borel-Cantelli, $\P(C_i \hbox{ i.o.})\allowbreak=1$. Note $$\frac{(t_{i+1}-t_i)\log\log\log (t_{i+1}-t_i)}{\log^2 (t_{i+1}-t_i)} =\frac{t_{i+1}\log\log\log t_{i+1}}{\log^2 t_{i+1}} \Big(1+o(1)\Big).$$ Also $$|V_i-R_{t_{i+1}}|+|\E V_i-\E R_{t_{i+1}}|\leq 2t_i=o\Big( \frac{t_{i+1}\log\log\log t_{i+1}}{\log^2 t_{i+1}}\Big).$$ Therefore with probability one, infinitely often  $$\ol R_{t_{i+1}} > \Big(1-\frac{\eta}{2}\Big)\xi t_{i+1}\log\log\log t_{i+1}/\log^2 t_{i+1}.$$ This proves the lower bound. 
\qed

We now turn to the  LIL for $-\ol R_n$. First we prove 
 
\begin{lemma}\label{LT5.2}  Let $\eps>0$. There exists $q_0(\eps)$ 
such that if 
$1<q<q_0(\eps)$, then 
$$\P(\max_{[q^{-1}n]\leq k\leq n} (\ol R_n-\ol R_k)>\eps n \log\log 
n/\log^2n) 
\leq \frac{1}{\log^2 n}$$ for $n$ large. 
\end{lemma} 
 
\proof Let $$G_k=(R_n-R_k)\frac{\log^2 n}{n}.$$ Let 
$$\sA_i=\left\{[q^{-1}n]+\Big[\frac{n\ell}{2^i}\Big]: \ell\in \Z_+\right\} 
\cap [0,n], \qquad i\leq \log_2 n+1.$$ Given $k$, let $k_i=\max\{j\in 
\sA_i: j\leq k\}$. We write 
$$G_k=G_{k_1}+(G_{k_2}-G_{k_1})+(G_{k_3}-G_{k_2})+\cdots,$$ where 
the sum is actually a finite one. If $\ol G_k>\eps \log\log n$ for some 
$[q^{-1}n]\leq k\leq n$, then either 
\begin{equation} \label{LE5.2} 
\ol G_{[q^{-1}n]} >\tfrac{\eps}{2} \log\log n 
\end{equation} 
  or for some 
$i$ there exist consecutive elements $\ell, m$ of $\sA_i$ such that 
\begin{equation}\label{LE5.3} 
\ol G_m-\ol G_\ell >\frac{\eps}{10i^2} \log\log n. 
\end{equation} 
 
By subadditivity $R_n-R_k\leq R_{n-k}\circ \theta_k$ for $k\leq n$, 
while by Lemma \ref{lemsecr21}
$$\E R_n-\E R_k\geq \E R_{n-k} -c_1 (1-q^{-1})^{1/2} \frac{n}{\log^2 
n}.$$ Then setting $k=[q^{-1}n]$, 
\begin{align*} 
\P(&\ol G_{[q^{-1}n]}>\tfrac{\eps}{2} \log\log n)\\ 
&=\P\Big((R_n-R_k)\frac{\log^2 n}{n}-(\E R_n-\E R_k)\frac{\log^2 n}{n} 
>\tfrac{\eps}{2} \log\log n\Big)\\ &\leq \P^{S_k}\Big( 
R_{n-k}\frac{\log^2 n}{n}-\E R_{n-k} \frac{\log^2 n}{n} 
+c_1(1-q^{-1})^{1/2} > \tfrac{\eps}{2}\log\log n\Big). 
\end{align*} 
  Using the fact that the $\P^x$ law of $R_{n-k}$ does not depend on $x$, 
this is the same as 
$$\P\Big(\frac{\ol R_{n-k}}{(n-k)/\log^2(n-k)}> \frac{n}{n-k} 
\frac{\log^2(n-k)} {\log^2 n}\Big( \tfrac{\eps}{2} \log \log 
n-c_1(1-q^{-1})^{1/2}\Big)\Big).$$ If $q>1$ is close enough to 1 and 
$n$ is large enough, by Theorem \ref{LU1.3} this is bounded by 
\begin{equation}\label{LE5.3a} 
  c_2 
\exp\Big(-c_3\frac{\eps}{2}\frac{1}{1-q^{-1}} \log \log n\Big)\leq 
\frac{1}{2\log^2 n}. 
\end{equation}  This bounds the probability of the event described in 
(\ref{LE5.2}). 
 
Similarly, $R_m-R_\ell \leq R_{m-\ell}\circ \theta_\ell$ and by Lemma 
\ref{lemsecr21}
$$\E R_m -\E R_\ell \geq \E 
R_{m-\ell}-c_1\Big(\frac{m-\ell}{n}\Big)^{1/2} 
\frac{n}{\log^2 n}.$$ So if $\ell$ and $m$ are consecutive elements of 
$\sA_i$,  similarly to (\ref{LE5.3a}) we obtain 
\begin{equation}\label{LE5.4} 
\P\Big(\ol G_m-\ol G_\ell\geq \frac{\eps}{10i^2}\log\log n\Big) 
\leq \P\Big(\frac{\ol R_{m-\ell}}{n/\log^2 n}\geq \frac{\eps}{10i^2} \log 
\log n-c_12^{-i/2}\Big). 
\end{equation}  For $n$ large, $c_12^{-i/2}\leq 
\frac{\eps}{20i^2}\log\log n$ for all $i$ and 
$n/(m-\ell)=2^i$, so by Theorem \ref{LU1.3} the left hand side of  (\ref{LE5.4}) is less than 
$$\P\Big(\frac{\ol R_{m-\ell}}{(m-\ell)/\log^2(m-\ell)} \geq 
\frac{\eps}{40i^2} \frac{n}{m-\ell}\log\log n\Big) 
\leq c_2 \exp\Big(-c_3 \frac{\log\log n}{40i^2}2^{i}\Big).$$ There are at 
most $2^{i+1}$ such pairs $\ell, m$, so 
\begin{align*} w_i&:=\P(\hbox{for some consecutive elements } \ell,m\in 
\sA_i: \ol G_m-\ol G_\ell >\frac{\eps}{10i^2}\log\log n) \\ &\leq c_2 
2^{i+1} \exp\Big(-c_3 \frac{\log\log n}{40i^2}2^{i}\Big). 
\end{align*} Since 
$c_32^i/40i^2 >2(i+1) \log 2$ for $i$ large, then for $n$ large enough 
$$w_i\leq c_2\exp\Big(-c_3 \frac{2^i\log\log n }{40i^2}\Big).$$ So then 
$$\sum_{i=1}^\infty w_i\leq \frac{1}{2\log^2 n}$$ for large $n$, and 
this bounds the event that for some $i$ there exist consecutive elements 
$\ell, m$ of $\sA_i$ such that (\ref{LE5.3}) holds. Combining with the 
bound for (\ref{LE5.2}), the  result follows. 
\qed

\noindent {\bf Proof of Theorem \ref{LT5.3}:} Let 
$$\Theta=(2\pi)^2 \det(\Gamma)^{-1/2} \kappa(2,2)^{-4}.$$ 
 
{\sl Upper bound.} Let $\eta, \eps>0$ and choose $q\in (1, q_0(\eps))$ 
where 
$q_0(\eps)$ is as in Lemma \ref{LT5.2}. Let $t_i=[q^i]$. If 
$$A_i=\left\{ -\ol R_{t_i}>(1+\eta)\Theta^{-1} \frac{t_i\log\log 
t_i}{\log^2 t_i}\right\},$$ then by Theorem \ref{LU1.3}, $\sum_i 
\P(A_i)<\infty$, and hence by Borel-Cantelli, 
$\P(A_i \hbox{ i.o.})=0.$ Let 
$$B_i=\left\{ \max_{t_i\leq k\leq t_{i+1}} (\ol R_{t_{i+1}}-\ol R_k) 
>\eps\frac{t_{i+1}\log\log t_{i+1}}{\log^2 t_{i+1}}\right\}.$$ By Lemma 
\ref{LT5.2}, $\sum_i \P(B_i)<\infty$, and again $\P(B_i\hbox{ i.o.})=0$. 
So with probability one, for $k$ large we have $t_i\leq k\leq t_{i+1}$ for 
some 
$i$ large, and then 
\begin{align*} -\ol R_k &=-\ol R_{t_{i+1}}+(\ol R_{t_{i+1}}-\ol R_k)\\ 
&\leq 
\Theta^{-1} (1+\eta) \frac{t_{i+1}\log\log t_{i+1}}{\log^2 t_{i+1}} + \eps 
\frac{t_{i+1}\log\log t_{i+1}}{\log^2 t_{i+1}}\\ &\leq 
q(\Theta^{-1}(1+2\eta)+2\eps) \frac{k\log\log k}{\log^2 k}. 
\end{align*}  Since 
$\eps, \eta$ can be made as small as we like and we can take $q$ as close 
to 1 as we like, this gives the upper bound. 
 
\medskip 
\noindent {\sl Lower bound.} Let $\eta>0$, 
$t_i=[\exp(i^{1+\frac{\eta}{2}})]$, 
$V_i=\#S((t_i,t_{i+1}])$. Let 
$$C_i=\left\{-\ol V_i\geq \Theta^{-1} (1-\eta) \frac{(t_{i+1}-t_i)\log\log 
(t_{i+1}-t_i)}{\log^2 (t_{i+1}-t_i)}\right\}.$$ By Theorem \ref{LU1.3}, 
$\sum_i 
\P(C_i)=\infty$. The $C_i$ are independent, and so by Borel-Cantelli, 
$\P(C_i\hbox{ i.o.})=1$. 
 
Since $R_{t_{i+1}}\leq V_i+R_{t_i}$ and $\E R_{t_{i+1}}\geq \E V_i$, 
then 
$$-\ol R_{t_{i+1}}\geq -\ol V_i-R_{t_i}.$$ Now 
$$R_{t_i}\leq t_i=o\Big(\frac{t_{i+1}\log\log t_{i+1}}{\log^2 
t_{i+1}}\Big)$$ and 
$$\frac{(t_{i+1}-t_i)\log\log(t_{i+1}-t_i)}{\log^2 (t_{i+1}-t_i)} 
\sim \frac{t_{i+1}\log\log t_{i+1}}{\log^2 t_{i+1}},$$ so 
$$-\ol R_{t_{i+1}}\geq \Theta^{-1}(1-2\eta) \frac{t_{i+1}\log\log 
t_{i+1}}{\log^2 t_{i+1}}, \qquad i.o.$$ This implies the lower bound. 
\qed

\def\noopsort#1{} \def\printfirst#1#2{#1} 
\def\singleletter#1{#1} 
     \def\switchargs#1#2{#2#1} 
\def\bibsameauth{\leavevmode\vrule height .1ex 
     depth 0pt width 2.3em\relax\,} 
\makeatletter 
\renewcommand{\@biblabel}[1]{\hfill#1.}\makeatother 
\newcommand{\bysame}{\leavevmode\hbox to3em{\hrulefill}\,}

\bigskip 
\noindent 
\begin{tabular}{lll}  & Richard Bass& Xia Chen\\ 
     & Department of Mathematics &  Department of Mathematics \\ 
    &University of Connecticut & University of  Tennessee\\ 
    &Storrs, CT 06269-3009 & Knoxville, TN 37996-1300 \\ 
&bass@math.uconn.edu&xchen@math.utk.edu\\ & &\\ & & 
\\ & & \\ 
     & Jay Rosen & \\ 
     & Department of Mathematics& \\ 
    &College of Staten Island, CUNY& \\ 
    &Staten Island, NY 10314& \\ &jrosen3@earthlink.net & 
\end{tabular}

\end{document}